 \documentclass[final,3p,times,twocolumn]{elsarticle}

\usepackage{amssymb}
\usepackage{amsmath}
\usepackage{xcolor}
\usepackage{algpseudocode}
\usepackage{algorithm2e}
\usepackage{float}
\usepackage{natbib}

\usepackage{gensymb}
\usepackage{textcomp}

\usepackage{subcaption}
\usepackage{bm}
\usepackage{booktabs}
\usepackage{multirow}

\journal{Ocean Engineering}

\begin{document}
\setlength{\parskip}{0pt}
\begin{frontmatter}



\title{MUTE‑DSS: A Digital-Twin-Based Decision Support System for Minimizing Underwater Radiated Noise in Ship Voyage Planning\\}

\author[1]{Akash Venkateshwaran}
\author[1]{Indu Kant Deo}
\author[1]{Rajeev K. Jaiman}


\affiliation[1]{organization={Department of Mechanical Engineering},
          addressline={The University of British Columbia},
           city={Vancouver},
           postcode={V6T 1Z4},
           state={BC},
           country={Canada}}



\begin{abstract}
We present a novel MUTE-DSS, a digital-twin-based decision support system for minimizing underwater radiated noise (URN) during ship voyage planning. It is a ROS2-centric framework that integrates state-of-the-art acoustic models combining a semi-empirical reference spectrum for near-field modeling with 3D ray tracing for propagation losses for far-field modeling, offering real-time computation of the ship noise signature, alongside a data-driven Southern resident killer whale distribution model. The proposed DSS performs a two-stage optimization pipeline: Batch Informed Trees for collision-free ship routing and a genetic algorithm for adaptive ship speed profiling under voyage constraints that minimizes cumulative URN exposure to marine mammals. The effectiveness of MUTE-DSS is demonstrated through case studies of ships operating between the Strait of Georgia and the Strait of Juan de Fuca, comparing optimized voyages against baseline trajectories derived from automatic identification system data. Results show substantial reductions in noise exposure level, up to 7.14 dB, corresponding to approximately an 80.68\% reduction in a simplified scenario, and an average 4.90 dB reduction, corresponding to approximately a 67.6\% reduction in a more realistic dynamic setting. These results illustrate the adaptability and practical utility of the proposed decision support system.


\end{abstract}




\end{frontmatter}


\section{Introduction}

The global shipping fleet has expanded significantly since the early 2000s and is projected to continue to grow rapidly due to several factors, including increased cargo volumes driven by the globalization of export industries in Asia, the adoption of larger and generally noisier vessels, and longer shipping routes. This expansion has contributed to a steady increase in ambient noise levels at low frequencies (10-100~Hz) across the ocean, with some studies reporting rates as high as 3~dB per decade \cite{mcdonald2006increases, andrew2011long}. By 2030, the maximum noise capacity of three primary segments of commercial shipping, namely container ships, bulk carriers, and oil tankers, is expected to increase by 87-- 102\%, driven by increased shipping activity, the use of larger vessels, and extended voyage distances \cite{kaplan2016coming}.

Several studies have shown that the impact of anthropogenic sound on marine organisms can vary widely, ranging from no effect at all to instant mortality, depending on the intensity, frequency and proximity of the noise to the source. Ship underwater radiated noise (URN), primarily peaking at low frequencies, significantly impacts mysticetes like whales that use similar sound ranges. This noise also extends to higher frequencies, harming odontocetes such as dolphins and porpoises. Masking is regarded as the most prevalent auditory impact of diffuse ship noise \cite{clark2009acoustic}. The exposure to ships' noise can induce a temporary or permanent reduction in auditory sensitivity, known as the noise-induced threshold shift, leading to hearing impairment \cite{southall2008marine}. Beyond this, ship noise can also cause various physiological and behavioral impacts in marine mammals, including elevated stress hormones, suppressed immune function, altered respiration rates, and changes in orientation \cite{Viola_Sciacca_2020,wright2007marine}. 

The urgency of addressing these impacts is exemplified by killer whales (\textit{Orcinus orca}) in British Columbia and Washington state, which have recently gained considerable attention regarding the effects of ship traffic due to their steadily declining population. Behavioral responses from killer whales have been observed at received noise levels exceeding 130~dB re 1~$\mu$Pa m (0.01--50~kHz), while the Lombard effect has been reported at ship noise levels above 98~dB re 1~$\mu$Pa m (1--40~kHz) \cite{williams2002behavioural, williams2014severity}. These studies indicate that killer whales may be more sensitive to acoustic disturbances than other species.


To address these challenges, several mitigation strategies have been developed that can be categorized into three groups: technological, operational, and economical \cite{vakili2020development}. Technological approaches encompass various noise-mitigation technologies addressing propeller, cavitation, and machinery noise, along with their respective levels of maturity and feasibility \cite{smith2022underwater}. Operational measures, such as voyage planning and optimal navigation strategies, can improve energy efficiency and reduce URN \cite{MEPC2023}. Voyage planning enables the use of alternative routes to avoid sensitive areas and allows for feasible speed reductions in designated conservation zones or during specific seasons, thereby mitigating URN impacts on marine ecosystems. The International Maritime Organization and the Convention for the Protection of the Marine Environment of the North-East Atlantic recognize speed reduction as an effective strategy for noise mitigation, which also significantly reduces the rate of collisions between ships and mammals \cite{MEPC2023, ospar2016ospar}. Exemplifying these efforts, the Vancouver Fraser Port Authority's Enhancing Cetacean Habitat and Observation (ECHO) Program has implemented voluntary vessel slowdowns in Haro Strait, Boundary Pass, and Swiftsure Bank, alongside re-routing ships in the Strait of Juan de Fuca \cite{ECHO_boundarypass, ECHO_swiftsure}. In particular, at Boundary Pass, these slowdowns resulted in a median noise reduction of 2.8 dB, a 48\% decrease in sound intensity and an approximate 18\% reduction in strike rates for both fin and humpback whales \cite{ECHO_boundarypass}.

The effectiveness of these operational strategies is difficult to quantify, as it depends on the type of vessel, the geographical environment, and the targeted species. In addition, such strategies can lead to increased fuel consumption, resulting in higher operational costs and emissions, necessitating a thorough trade-off analysis before implementing any such measures. Hence, there is a critical need for a management system that considers these factors for effective decision-making. This paper presents a decision support system for ship voyage optimization that addresses the complex trade-offs between operational efficiency, economic considerations, and marine mammal conservation, providing a comprehensive framework for mitigating ship noise impacts on vulnerable marine ecosystems.

Recent research covers a variety of methods, including dynamic programming, genetic algorithms, and machine learning, to address speed optimization, weather routing, and operational decision making. Dynamic programming is widely used to optimize both route and speed, considering weather forecasts, fuel consumption, and ship comfort by discretizing the voyage into stages and using weather data to minimize fuel use and improve safety \cite{zaccone2018ship}. Graph-based approaches, such as Dijkstra's algorithm and its 3D variants, enable multi-objective optimization (\textit{e.g.}, minimizing fuel, damage, and ensuring accurate arrival times) \cite{bottner2007weather,wang2019three}. Artificial neural networks and integrated simulation models are used to predict fuel consumption under various operational conditions, allowing for real-time adjustments in speed and trim to maximize efficiency \cite{sang2023ship,du2019two}. Combining genetic algorithms with dynamic programming allows for more flexible and robust optimization, especially in dynamic sea environments, leading to significant reductions in fuel consumption and emissions \cite{wang2019three,wang2021voyage}. Some studies incorporate uncertainty in weather forecasts into the optimization process, using risk penalty functions to avoid underestimated hazards and improve both safety and energy efficiency \cite{chen2024isochrone,yuan2022uncertainty}. The choice of energy cost models (empirical, data-driven, or hybrid) significantly affects optimization outcomes, making accurate modeling of the speed-energy relationship crucial for reliable voyage planning \cite{mylonopoulos2023comprehensive}.

The existing literature on ship voyage planning that considers URN emissions remains limited, with only two significant contributions identified. The first study formulates the problem with URN constraints specifically for all-electric ships, focusing on optimal scheduling of generators and energy storage devices to maintain desired speeds while satisfying URN limits and ensuring reliable operation of the ship's electric power network \cite{khatami2023optimal}. However, this work employs a simplified empirical source spectrum model, Ross's model \cite{Ross1976}, which determines noise levels through a power-law relationship with vessel speed, potentially unrepresentative of modern vessel configurations. The second contribution presents a multi-objective optimization framework for fixed-route vessel speed optimization that integrates ship noise signature modeling with acoustic propagation models to generate 2D noise maps \cite{venkateshwaran2024multi}. While this study demonstrates a successful reduction in overall noise intensity levels without significantly impacting fuel consumption, several critical limitations restrict its practical application. The optimization process assumes a 2D ocean environment with stationary marine mammal locations, failing to account for the dynamic nature of marine ecosystems. Additionally, the usage of Ross's noise spectral model make the findings less representative of contemporary vessel noise characteristics.

This paper addresses these limitations through the development of the Mitigating Underwater Noise Transmission and Effects Decision Support System (MUTE-DSS). The proposed system reformulates ship voyage planning as a two-stage optimization process that encompasses both route planning and speed profile optimization while incorporating considerations of environmental impact. The MUTE-DSS framework enables 3D simulation of ship voyage operations in scenarios where vessels transit between ports (e.g., Port A to Port B) while encountering multiple marine mammals along their route. Key innovations include: (1) integration of data-informed, dynamically moving mammal populations within the simulation environment; (2) implementation of real-time acoustic ship noise mapping models; and (3) replacement of the outdated Ross model with contemporary, robust ship noise signature models that accurately represent modern vessel characteristics. The optimization process quantifies the impact of ship noise through objective functions based on the received noise levels experienced by individual marine mammals within the simulation environment. This approach enables adaptive adjustment of both route and vessel speed profiles to minimize the impact of URN with voyage constraints.  Moreover, MUTE-DSS is built on a ROS2-centric architecture \cite{macenskirobot} that interfaces directly with external control systems and sensors, paving the way for seamless hardware-in-the-loop integration.

This ROS 2 foundation is pivotal in positioning MUTE-DSS as a maritime digital twin (DT): an ecosystem in which a live, bidirectional data loop synchronizes a physical vessel with its virtual counterpart in real-time. DT capabilities such as real-time data exchange, predictive analytics, diagnostics, and feedback extend classical simulation beyond static physics and graphics engines \cite{semeraro2021digital}. The ROS2 publish-subscribe middleware, standardized message types, and extensive driver ecosystem ensure temporally coherent data flows between sensors, controllers, and virtual processes across multiple robots or vessels. Within MUTE-DSS, ROS nodes stream sensor topics (e.g., IMU, sonar, DVL) into the twin, while the twin’s analytics publish commands or parameter updates back to the physical platform, closing the loop. Thus, by having DT principles via ROS2, MUTE-DSS transforms voyage planning into a predictive self-optimizing process that mitigates underwater-noise pollution in dynamically evolving marine ecosystems.

The problem description and its underlying assumptions are detailed in Section \ref{sec_math_prob}. In Section \ref{sec_math_forml}, we introduce the datasets employed in this study and formally define three principal unit of the MUTE-DSS framework: First, the Modeling Unit (Section \ref{sec_modeling_unit}) comprises two subcomponents: ship noise characterization and marine mammal modeling. Second, the Optimization Unit (Section \ref{sec_optimization_unit}) comprises route planning and speed optimization. Finally, the Simulation and Visualization Unit (Section \ref{sec_simulation_visualization_unit}) integrates the output of the preceding units to generate comprehensive scenario analyses and graphical representations of ship-mammal interactions. Section \ref{sec_case_studies} presents the case studies, and their evaluation is reported in Section \ref{sec_results}. Finally, Section \ref{sec_conclusion} offers concluding remarks and outlines directions for future research.

\begin{figure*}[ht]
  \begin{center}
    \includegraphics[width=0.75\textwidth]{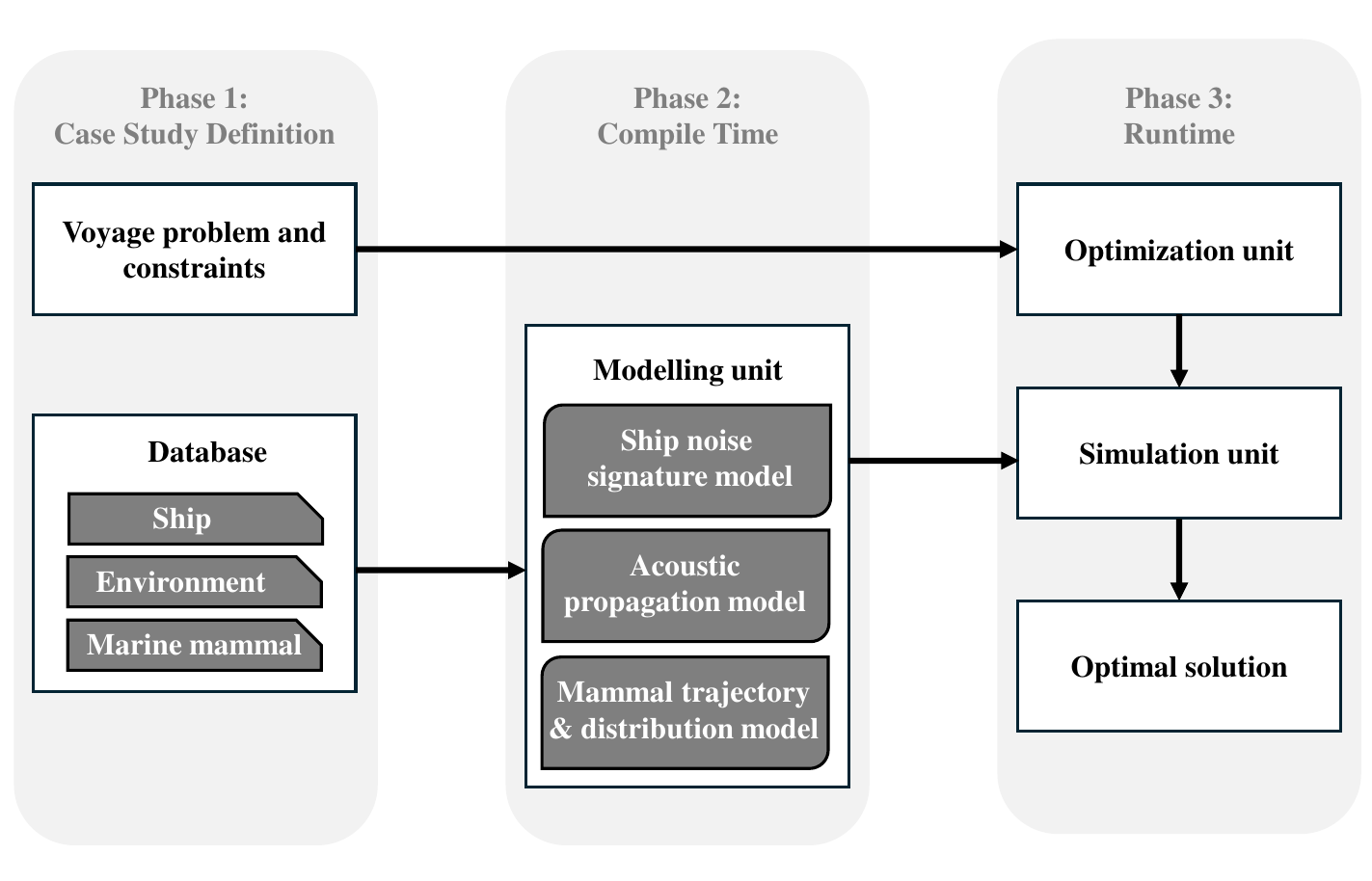}
    \caption{
     Schematic overview of the MUTE-DSS for ship voyage optimization. Three phases of ship voyage optimization are shown: case study definition (Phase 1), model compilation (Phase 2), and dynamic runtime simulation and optimization (Phase 3)
    }
    \label{fig_flowchart}
  \end{center}
\end{figure*}

\section{Problem description}\label{sec_math_prob}
While numerous studies have addressed individual aspects of voyage optimization and URN modeling, few offer an integrated system that accounts for realistic ship operations and ecological impact in a dynamically evolving environment. To this end, we now define the problem formulation and overall architecture of the MUTE-DSS framework to address the impact of URN from ships. The problem is formulated as a single-ship voyage optimization problem where URN levels serve as the primary objective function. Optimization integrates route planning and speed control to reduce the overall acoustic impact on marine ecosystems within a simulated 3D environment.
The ship operates within a spatial region that contains moving marine mammals, dynamically adjusting speed and trajectory based on real-time acoustic conditions. This adaptive approach requires comprehensive models that capture the acoustic impact of various voyage scenarios and their interactions with marine life. The modeling framework comprises three sequential phases, as illustrated in Fig. \ref{fig_flowchart}.

Phase 1 establishes the voyage problem by specifying permitted shipping lanes and enforcing exclusion zones (e.g., coastlines, shallow waters, ecologically sensitive areas). Operational constraints, such as maximum speed limits and target arrival time, are also defined to reflect realistic voyage conditions. Once the voyage parameters and constraints are in place, three categories of data are utilized: (1) ship characteristics (e.g., engine specifications, hull geometry, dimensions) for URN modeling; (2) environmental inputs (e.g., bathymetry, temperature, salinity profiles) to model URN propagation; and (3) marine mammal information (e.g., species distributions, behavioral patterns, movement trajectories) for impact assessment.

In Phase 2, the data from Phase 1 are used to develop computational models that support online optimization and simulation. A source-level model predicts ship noise spectra based on vessel characteristics. An acoustic propagation model then rapidly estimates three-dimensional transmission loss across the study area. Finally, a marine mammal trajectory and distribution model simulates movement patterns within the 3D environment. Together, these components enable real-time evaluation of ship–mammal interactions during voyage simulations. Further specifications are provided in the following section.

Phase 3 integrates all models into a dynamic framework to simulate a realistic ship voyage scenario. This phase includes an optimization unit that performs route planning and speed optimization, guiding the vessel to optimal operating conditions. Following this, the simulation unit reconstructs the optimized voyage and interaction with the marine mammals, allowing the solution to be interpreted and visualized.

The MUTE-DSS relies on several key assumptions, which are stated as follows:
\begin{enumerate}
    \item \textbf{Weather conditions}: Environmental factors such as wind and waves introduce additional resistance during navigation, reducing vessel speed compared to calm water conditions and impacting fuel consumption and route planning \cite{vitali2020coupling}. The current MUTE-DSS focuses on ship-radiated noise mitigation strategies; therefore, these effects are not active, but the architecture is designed to accommodate them in future extensions via modular, physics-based disturbance models (e.g., wind, wave, and current generators). It can integrate water surface interaction simulators (e.g., WaveSim \cite{santos2006wavesim} or Kelpie \cite{mendoncca2013kelpie}) through its client interface to supply context-aware dynamics and environmental forces such as water resistance and buoyancy.

    \item \textbf{Uncertainty considerations}: Maritime voyage planning involves significant uncertainties, with only 55-89\% of vessels arriving as scheduled \cite{elmi2022uncertainties}. Port arrival and transit times are affected by traffic congestion, cargo handling delays, and adverse weather conditions \cite{notteboom2006time}. While uncertainty quantification is not addressed in this study, it represents a potential area for framework enhancement.
    
    \item \textbf{Problem scope}: The primary objective function focuses on the impact of URN. Additional practical navigation considerations-including vessel motions, collision avoidance, navigation risk, and maritime service fees-are not modeled but could be incorporated as supplementary objectives or constraints. Kinematic and dynamic constraints associated with ship motion and maneuvering are similarly excluded. The MUTE-DSS framework provides sufficient flexibility to accommodate these factors in future implementations.
    
    \item \textbf{Model assumptions}: Ship noise signature modeling and acoustic propagation modeling assumptions are detailed in Sections~\ref{sec_acoustic_modeling}. Additionally, the mammal trajectory model is simplified, with each marine mammal assigned a random linear velocity vector, providing an initial approximation of movement.
    
    \item \textbf{Temporal variations}: To enable computation of real-time objective functions, noise modeling is precomputed as described in Phase 2 of Figure~\ref{fig_flowchart}. Consequently, the environmental parameters of the acoustic model remain constant throughout each voyage, and seasonal variations of oceanographic conditions are not considered in the current implementation.
\end{enumerate}

\section{Mathematical modeling}\label{sec_math_forml}
\begin{figure*}[!ht]
  \centering
  \makebox[\textwidth][c]{%
    \includegraphics[width=0.9\textwidth]{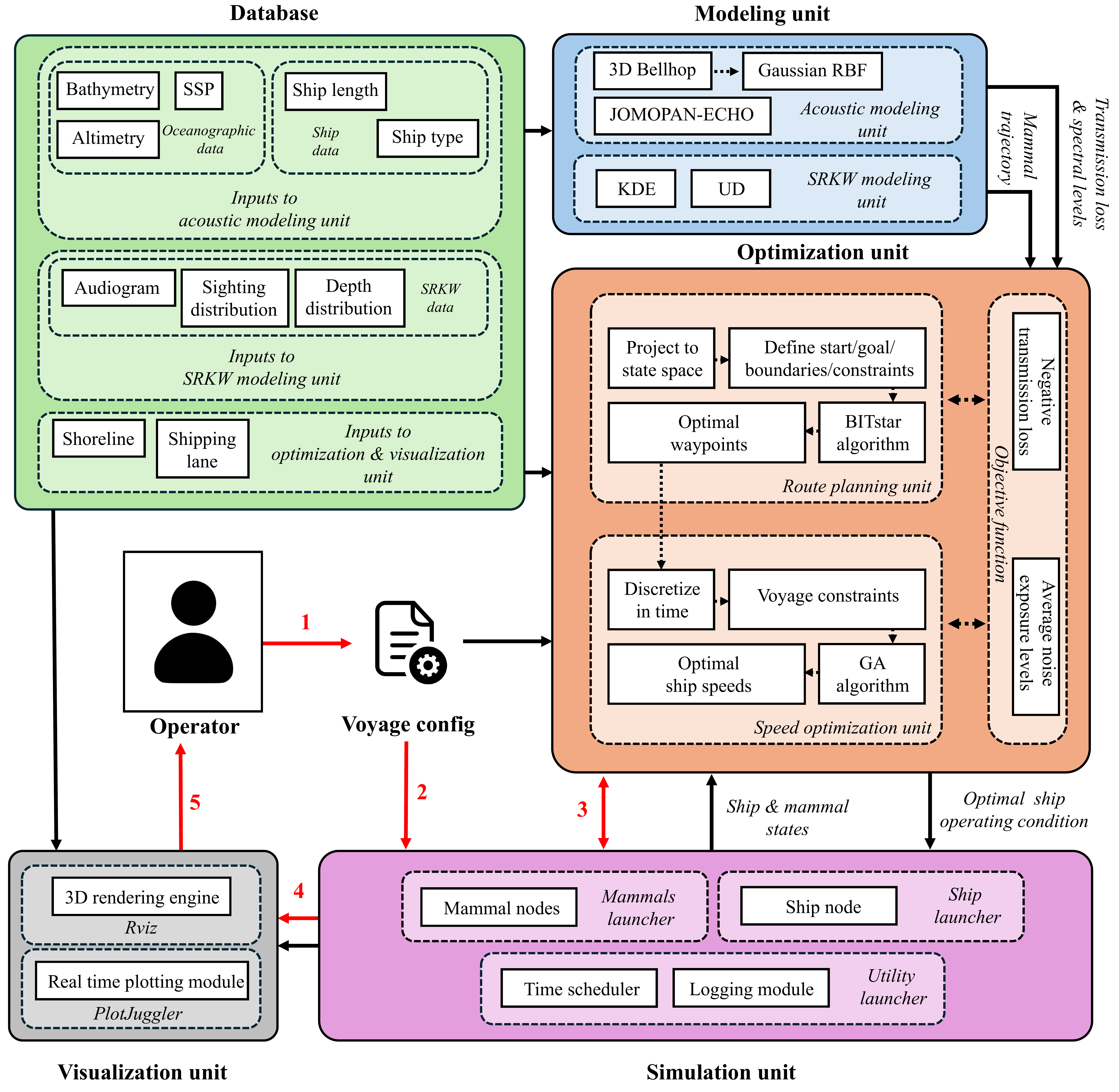}}
\caption{Digital twin architecture of the proposed MUTE-DSS framework, showing interaction among modeling, optimization, simulation, and visualization units. Red arrows indicate system data flow from voyage initialization to real-time decision support. }

  \label{fig_pipeline}
\end{figure*}

This section presents the underlying models, objective functions, and constraints relevant to the optimization problem. Fig. \ref{fig_pipeline} illustrates the architecture of the MUTE-DSS, with comprehensive explanations provided for each component's mathematical formulation.

The workflow, indicated by the red arrow in Fig. \ref{fig_pipeline}, starts with the operator providing the configuration of the trip: (1) This configuration file encompasses details regarding the voyage, such as the departure coordinate, destination coordinate, and estimated time of arrival (ETA), along with details related to marine mammals, i.e. the number of mammals distributed within the 3D environment. Subsequently, (2) the simulation unit initializes the ship as a node at the designated departure coordinate, disperses the specified number of mammal nodes per the modeling unit, and activates the time scheduler and logging units. Following this, (3) the optimization unit, which identifies both the ship and mammal nodes, is responsible for guiding the ship's optimal operating conditions towards the destination. This unit integrates both the speed optimization and route planning units. It retrieves the states of the mammals and ship nodes, computes the voyage cost based on the objective function defined by the modeling unit, and optimizes the trajectory over a specified time window or until convergence criteria are met. The simulation unit, in turn, receives feedback on the optimal ship's speed and route with respect to the simulation time, forming a continuous feedback loop. Finally, (4) the operation scenario involving the ship and mammals in the 3D environment is visualized through the interactive visualization unit. In conclusion, (5) the logged data and simulated operating conditions for the specific voyage configuration are provided to the operator.
The individual components of the architecture are explained in the following subsections.

\subsection{Database}\label{sec_database}
The database unit, as shown in Fig. \ref{fig_pipeline}, is categorized into data necessary for acoustic modeling and data pertinent to the modeling of marine mammals. The study area encompasses the waters of the Juan de Fuca Strait to the Strait of Georgia. Data collection and simulated case studies are restricted to this specified region. This region overlaps with key international shipping routes connecting the Pacific Ocean to major ports in southern British Columbia, including the Port of Vancouver, the largest in the region. In addition, it intersects with the Washington State Ferries route that runs from Anacortes, Washington, to Sidney, British Columbia.

\begin{figure}[ht]
\centering
\includegraphics[width=0.45\textwidth]{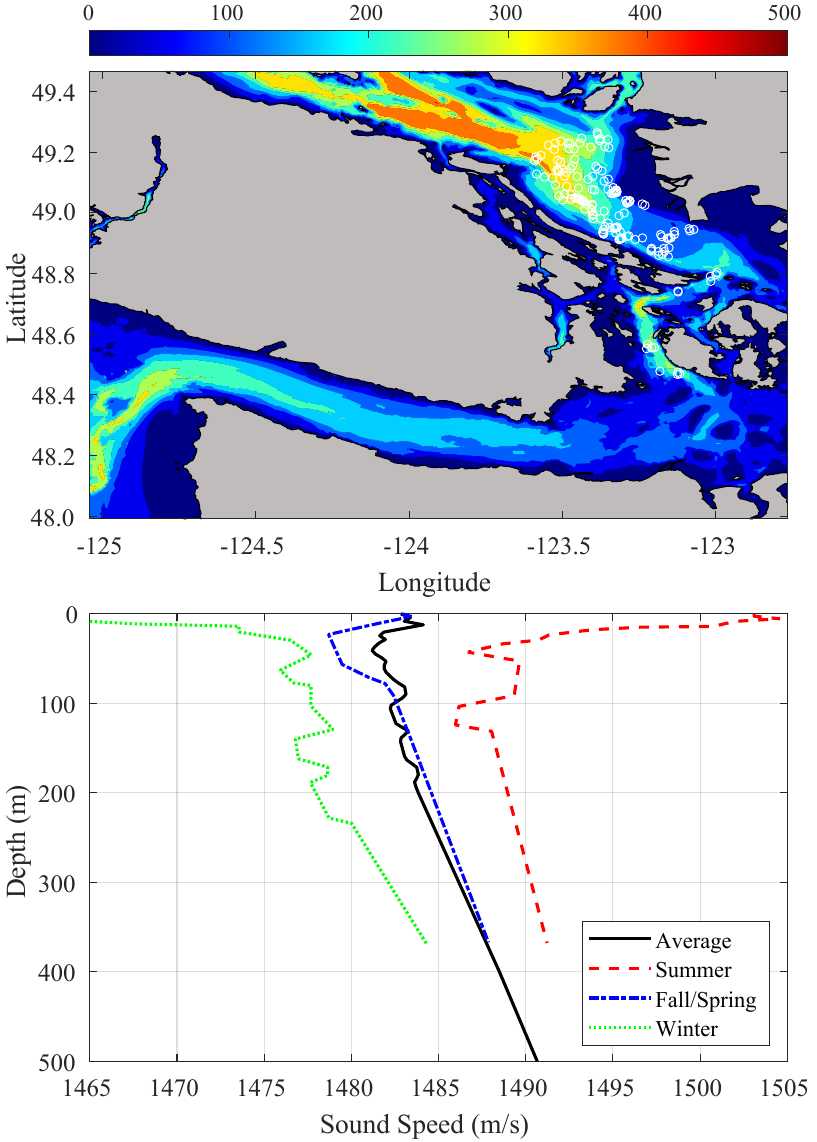}
\caption{Acoustic database for the study area: (a) bathymetry contours with SSP measurement locations (white circles); (b) seasonal SSPs at the site from \cite{Eickmeier2021SalishSea} and the average SSP used in this study.}
\label{fig_database}
\end{figure}

Key factors such as bathymetry, surface roughness, and seabed composition (e.g., sand, silt, and clay) are crucial in sound propagation. Moreover, sound propagation is significantly influenced by the refractive properties of SSP as they vary with depth. The bathymetry data for the study region have been sourced from the GEBCO database \cite{gebco2024}. The seasonal variations of SSP in the Salish Sea are obtained from \cite{Eickmeier2021SalishSea}, providing seasonal measured profiles in the Strait of Georgia. Fig. \ref{fig_database} presents the bathymetry contour of the study area and the average SSP used for the analysis.

The marine mammal modeling requires a density model to initialize the mammals within the environment, as well as a trajectory model to estimate their positions over time after initialization. In this study, Southern resident killer whales (SRKWs) are considered the sole target species, and the necessary models are specifically designed for them. SRKW is classified as "Threatened" under both the Canadian Species at Risk Act and the U.S. Endangered Species Act. SRKWs inhabit the Salish Sea throughout the year, with a significant concentration in Haro Strait, which lies off the west coast of San Juan Island and forms the core of their critical habitat \cite{thornton2022southern}. More details on the dataset and models are presented in Section \ref{sec_srkw_modeling}. 
The database described above provides the foundation for environmental and species-specific modeling. We now describe how these data sets are leveraged in the modeling unit to characterize ship noise emissions, their underwater propagation, and the movement of marine mammals.

\subsection{Modeling unit}\label{sec_modeling_unit}
This section delineates the modeling unit, which uses the database as input. As previously outlined, model development is conducted prior to runtime, specifically during the compilation phase, as illustrated in Fig. \ref{fig_flowchart}. The modeling unit itself comprises two subcomponents: acoustic modeling and SRKW behavior modeling.

\subsubsection{Acoustic modeling}\label{sec_acoustic_modeling}

This section addresses the modeling of URN from ships. The noise level at the source (NLS) or near-field noise levels of the ship is the sound pressure level (SPL) measured at a reference distance of one meter from the source. It is independent of the environment in which the source operates and is quantified as the logarithmic ratio of the acoustic intensity produced by the source to a standard reference intensity. Meanwhile, the noise level received by mammals (NL), or far-field noise levels, is lower than the NLS due to transmission loss. It can be determined using a passive sonar equation, as follows:
\begin{equation}
\text{NL} = \text{NLS} - \text{TL},
\end{equation}
where all three quantities are in dB ref 1 $\mu$Pa. Transmission loss (TL) quantifies the reduction in the acoustic intensity at a given location relative to the source intensity and represents the dissipation of energy influenced by several factors, including the loss attributed to geometric spreading, boundary effects, scattering phenomena, and volume attenuation. 

Ship noise comprises broadband and tonal components generated by three main sources: machinery (including propulsion and auxiliary equipment), hydrodynamic flow interactions with the hull, and propeller operations \cite{mckenna2012underwater}. To model these sources, many empirical and semi-empirical ship noise models have been developed since WWII, ranging from Ross's foundational 1976 model \cite{Ross1976} to the recent physics-based PIANO model \cite{lloyd2024piano}. These models incorporate various ship parameters such as speed, displacement, length, and propeller characteristics to predict underwater radiated noise levels. Their comparative analysis is presented in \ref{sec_app_acoustic}.

In this study, we have implemented the JOMOPANS-ECHO (JE) model, which is based on regression analysis of monopole ship source level measurements from the ECHO Program \cite{jompans_echo2021}. A total of 1,862 source level measurements were collected by JASCO Applied Sciences in the Haro Strait and Strait of Georgia during the 2017 ECHO slowdown trial \cite{macgillivray2019slowing}. The model incorporated ship-type-specific reference speeds and spectrum coefficients, as the ship type is a key factor in determining URN, as mentioned in the previous section. All necessary information about the ships was recorded using the automated identification system (AIS). The source levels from the hydrophone measurements were calculated using a frequency-dependent propagation loss model, rather than the spherical propagation model used in previous source models.

The JE model reintroduced the relationship with ship length and proposed the following equation:
\begin{equation}
\begin{aligned}
\mathrm{NLS}_{JE}(f,v,l,T) &= \mathrm{NLS}_{0}(\hat{f},T) \\
 & + 60 \log \left(\frac{v}{v_{T}}\right) + 20 \log \left(\frac{l}{l_{0}}\right),
\end{aligned}
\label{eq_jomopan}
\end{equation}
where $v$ and $l$ are the ship speed in knots and length in feet, respectively. The reference speed ($v_{T}$) based on the ship type and the reference length ($l_0$) has a value of 300 ft. The baseline spectrum ($\mathrm{NLS}_{0}(\hat{f},T)$) is defined as:
\begin{equation}
\begin{aligned}
    \mathrm{NLS}_{0}(\hat{f},T) &= K - 20 \log(\hat{f}_1) \\
    & - 10 \log \left( 1 - \left( \frac{\hat{f}}{\hat{f}_1} \right)^2 + D^2 \right),
\end{aligned}
\end{equation}
where \( \hat{f} = \frac{f}{f_{ref}}, \hat{f}_1 = 480 \, \text{Hz} \times \left( \frac{v_{ref}}{v_{T}} \right), f_{ref} = 1 \, \text{Hz}, v_{ref} = 1 \, \text{kt}, K = 191 \, \text{dB}, D = 3 \) for all the ship types except the cruise vessel having $D = 4$.
For cargo vessels (container ships, vehicle carriers, bulkers, tankers), the baseline spectrum below 100 Hz is changed to:
\begin{equation}
\begin{aligned}
\mathrm{NLS}_{0}(\hat{f}, T) &= K - 40 \log \left( \frac{\hat{f}}{\hat{f_1}} \right) + 10 \log(\hat{f}) \\
& - 10 \log \left( \left(1 - \left( \frac{\hat{f}}{\hat{f_1}}\right)^2\right)^2 + D^2 \right) \, \text{dB}
\end{aligned}
\end{equation}
where \( \hat{f}_1 = 600 \, \text{Hz} \times \left( \frac{v_{ref}}{v_{T}} \right), K = 208 \, \text{dB}, D = 0.8 \) for container ship and bulker or $D = 1.0$ for vehicle carriers and tankers.

 The JE model includes an estimated statistical uncertainty of ±6 dB (rms) across a 0.02–20 kHz frequency range. Additionally, the model explicitly accounts for the speed-dependent nature of source levels, following a well-supported power-law trend for cavitation noise, which is critical for evaluating slow-down mitigation strategies. Moreover, the model parameters are tuned based on measurements in the study area where this analysis is conducted, and they account for various ship types in the Salish Sea, making them an ideal candidate for noise source modeling. Further details of the comparison noise source models and reference table of the ship types for the JE model are presented in  \ref{sec_app_acoustic}.

 The modeling of underwater acoustic propagation involves solving linearized hydrodynamic equations that account for sound speed profiles (SSPs), pressure, particle velocity, and density. A propagation model is fundamentally a computational tool used to calculate the pressure field by solving these equations for a specified set of transmitted signals, propagation speed, surface and bottom characteristics, among other factors. There are two broad classes of solutions: numerical models that use techniques like parabolic equations, normal modes, and ray solutions, etc \cite{wang2014review}; and data-driven models such as convolutional recurrent autoencoder networks \cite{mallik2022predicting} and range-dependent conditional convolutional neural networks \cite{deo2025predicting} that employ deep learning techniques.

\begin{figure*}[!ht]
  \begin{center}
    \includegraphics[width=\textwidth]{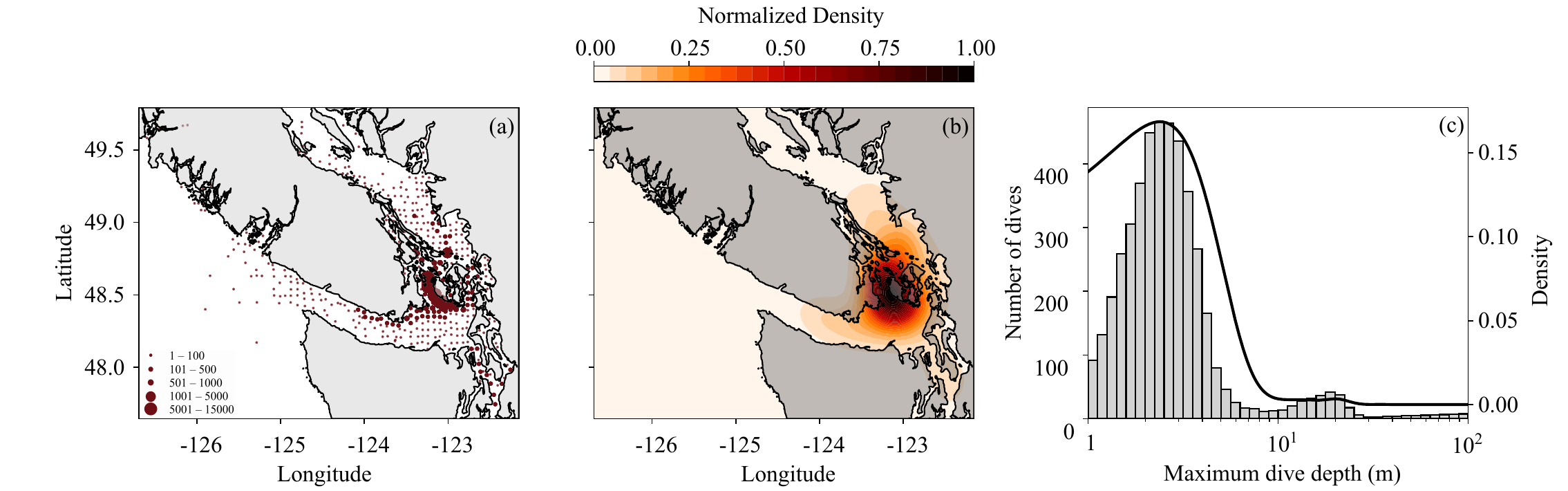}
    \caption{SRKW database of the study area: (a) sightings and encounter distribution with point size proportional to sighting per 5 km$^2$ \cite{thornton2022southern}; (b) kernel density fit of sighting data; (c) histogram of maximum dive depths \cite{mcrae2024killer}.}

    \label{fig_SRKW}
  \end{center}
  
\end{figure*}

In MUTE-DSS, we used 3D Bellhop \cite{acoustics_toolbox}, a highly efficient ray tracing program that is based on the theory of Gaussian beams, to compute TL. The NLS and TL are computed at the central frequency for each third-octave band within the 12.5 Hz to \(10^4\) Hz range. Given that the JE model is not computationally intensive, the source noise spectrum is calculated dynamically during the optimization process. Conversely, for TL, the 3D Bellhop model is precomputed for each specified frequency using UBC Sockeye HPC clusters. A total of 354 source locations are strategically positioned along the ship lanes of the study area, with TL calculated within a 100 km radius from each source location. The depth analyzed extends up to 100 meters, aligning with the distribution of SRKWs, which predominantly occupy this depth range, as depicted in Fig. \ref{fig_SRKW}. Under the simplifying assumption that the ship acts as a point-like source of sound, a single source depth of 6 m is considered in this study.

\begin{figure*}[ht]
  \begin{center}
    \includegraphics[width=0.8\textwidth]{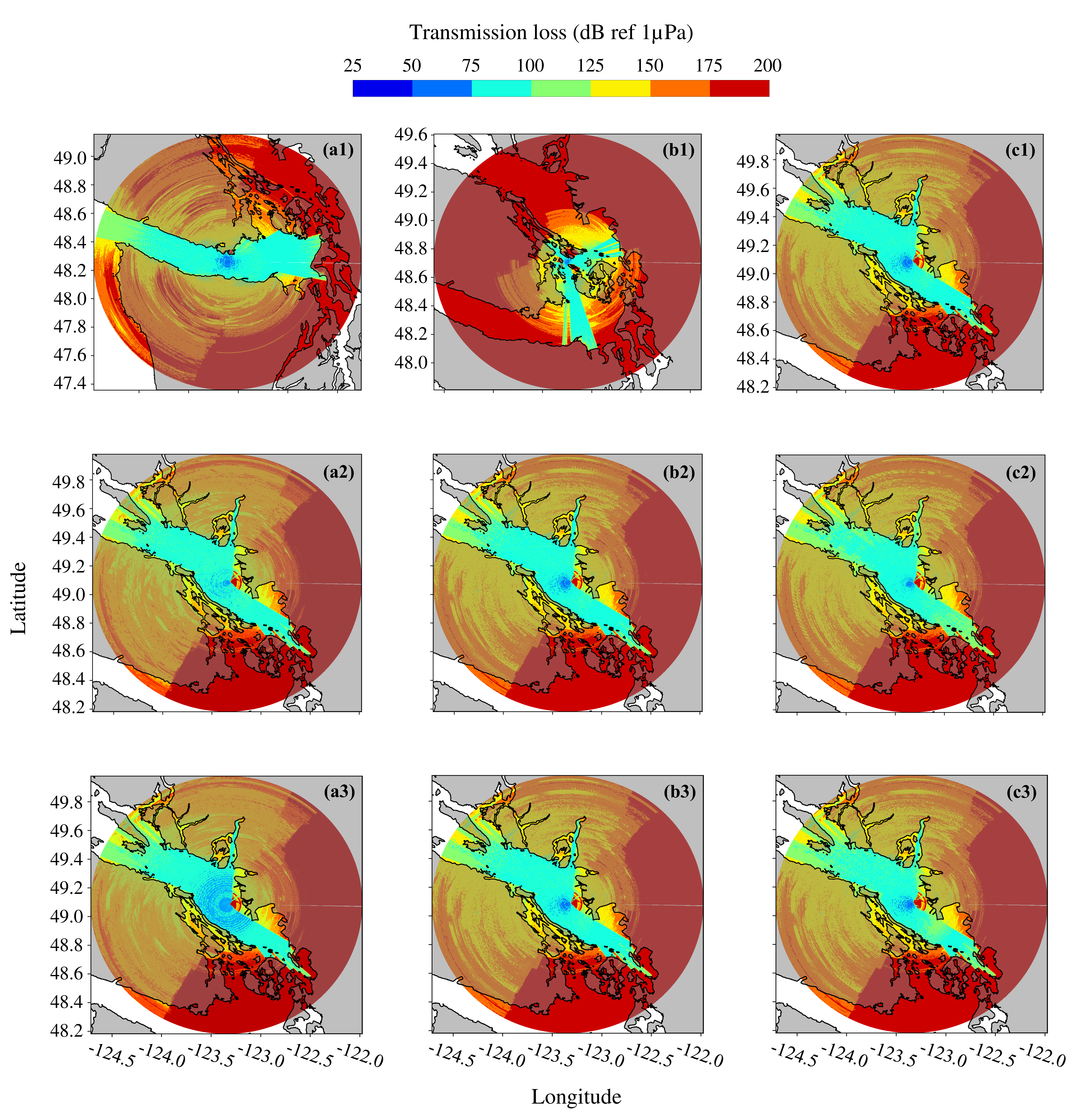}
    \caption{Transmission loss contour map for the study area, illustrating variations across multiple parameters: (1) source locations; (2) frequencies of 12.5 Hz, 100 Hz, and 1000 Hz; (3) depth sections at 0 m, 50 m, and 100 m.}
    \label{fig_transmission_loss}
  \end{center}
  
\end{figure*}

Figure \ref{fig_transmission_loss} presents the TL contours for three configurations: different source locations in the first row, varied frequencies in the second row, and different depth planes (z-planes) in the third row. We require a fast function approximation to compute TL based on the source location (ship node) and receiver location (mammal nodes). This multivariate interpolation problem is addressed by using radial basis function (RBF) interpolation. 
RBF interpolation, originally introduced by Hardy in 1971 \cite{hardy1971multiquadric}, is a widely used technique for multivariate interpolation, particularly effective in handling scattered data interpolation problems. In this study, a Gaussian RBF is used and is found to be stable for the TL approximation function. 

Let the data set be $(\mathbf{x}_i,\mathrm{TL}_i)$ for $i=1,2,\dots,n$, where each
\[
\mathbf{x}_i
  = \bigl(s_{\mathrm{lat}},\,s_{\mathrm{long}},\,r_{\mathrm{lat}},\,r_{\mathrm{long}},\,r_z\bigr)^{\top}
  \in\mathbb{R}^{5},
\]
where $s_{\mathrm{lat}}$ and $s_{\mathrm{long}}$ represent the source latitude and longitude, while $r_{\mathrm{lat}}$, $r_{\mathrm{long}}$, and $r_z$ represent the receiver latitude, longitude, and depth, respectively. An interpolant $\widetilde{\mathrm{TL}}\colon\mathbb{R}^{5}\!\to\!\mathbb{R}$ that satisfies
\begin{equation}
   \widetilde{\mathrm{TL}}(\mathbf{x}_i)=\mathrm{TL}_i, \qquad i=1,2,\dots,n.
\label{eq_interpolant_cond}
\end{equation}
In the case of RBF interpolation, we require that $\widetilde{\mathrm{TL}}(\mathbf{x})$ be a weighted linear combination of kernel functions $\psi(r)$, i.e.

\begin{align}
    \widetilde{\mathrm{TL}}(\mathbf{x}) &= \sum_{i=1}^n \lambda_i \psi(r_i), \quad \psi(r_i) = \exp\left(-0.5\left(\frac{r_i}{\sigma}\right)^2\right),
\label{eq_interpolant}
\end{align}
where $r_i = \|\mathbf{x} - \mathbf{x}_i\|$ is the Euclidean distance between $\mathbf{x}$ and $\mathbf{x}_i$, and $\sigma$ is a parameter that controls the width of the Gaussian function. 
From Eq.~\ref{eq_interpolant_cond} and Eq.~\ref{eq_interpolant}, it follows that
\begin{equation}
\begin{aligned}
\sum_{i=1}^{n} \lambda_i \psi(\| \mathbf{x}_j - \mathbf{x}_i \|) &= \mathrm{TL}_j,  \\
 j &= 1, 2, \dots, n, \\
\begin{bmatrix}
\psi(\|\mathbf{x}_1 - \mathbf{x}_1\|)  & \cdots & \psi(\|\mathbf{x}_1 - \mathbf{x}_n\|) \\
\psi(\|\mathbf{x}_2 - \mathbf{x}_1\|) & \cdots & \psi(\|\mathbf{x}_2 - \mathbf{x}_n\|) \\
\vdots & \ddots & \vdots \\
\psi(\|\mathbf{x}_n - \mathbf{x}_1\|) & \cdots & \psi(\|\mathbf{x}_n - \mathbf{x}_n\|)
\end{bmatrix}
\begin{bmatrix}
\lambda_1 \\
\lambda_2 \\
\vdots \\
\lambda_n
\end{bmatrix} &= 
\begin{bmatrix}
\mathrm{TL}_1 \\
\mathrm{TL}_2 \\
\vdots \\
\mathrm{TL}_n
\end{bmatrix}, \\
\boldsymbol{\Psi} \boldsymbol{\lambda} &= \boldsymbol{\mathrm{TL}}.
\end{aligned}
\label{eq_interpolant_sys}
\end{equation}

\begin{figure*}[ht]
  \begin{center}
    \includegraphics[width=\textwidth]{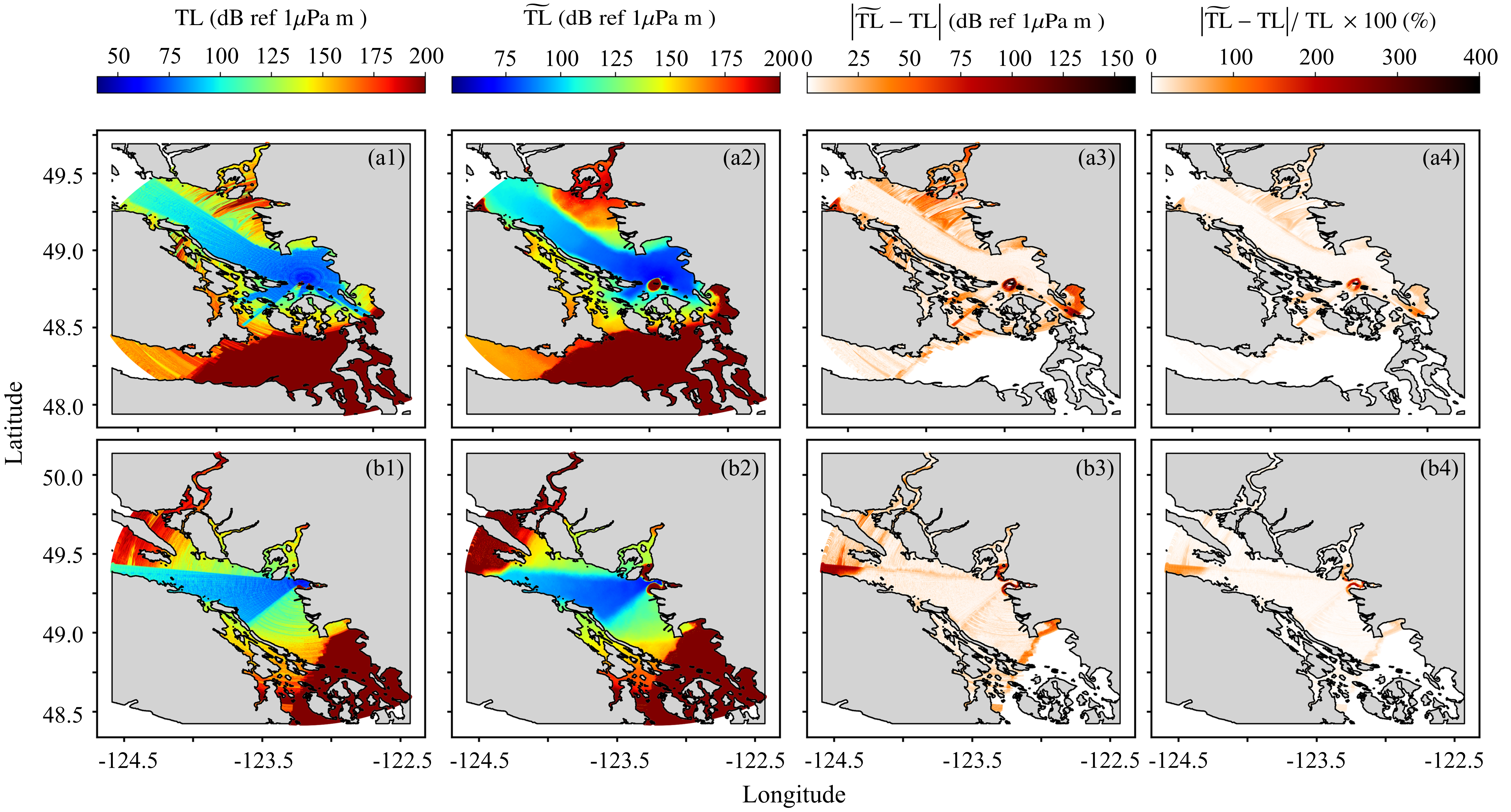}
    \caption{Comparison of transmission loss data and model prediction at 12.5 Hz: (1) 3D Bellhop data computed over a 100 km radius; (2) prediction from the Gaussian RBF interpolant; (3) absolute error; (4) absolute percentage error (\%).}
    \label{fig_RBF}
  \end{center}
\end{figure*}

Then, this system of equations is solved to determine the weights ($\boldsymbol{\lambda}$) of the approximation function. 
The squared interpolation matrix ($\mathbf{\Psi}$) in Eq. \ref{eq_interpolant_sys} has dimensions equal to the number of data points, making it computationally infeasible to consider the entire 3D dataset of varying source and receiver locations in a naive manner. To mitigate this, data are initially clustered and only a limited subset is sampled from each cluster, effectively reducing the dimensions of the matrix without compromising accuracy. Several trials are conducted to fine-tune the clustering and sampling process to get to the final interpolant. Each component of the TL represents a one-third octave band for each frequency, covering a total of 30 components. These 30 components are reduced to 10 using principal component analysis, which significantly improves the model's speed in computing TLs across frequencies. Subsequently, each of these 10 components is modeled using a separate RBF interpolant. Fig. \ref{fig_RBF} compares the Bellhop TL contours, Gaussian RBF interpolant predictions, absolute error, and percentage absolute error. While the contour aligns with the data, some artifacts emerge, particularly in regions of large TL variation. These artifacts are exaggerated and spatially dispersed due to the fixed width of the Gaussian kernel. 

The JOMOPANS-ECHO model and the Gaussian RBF interpolant described in this section enable fast computation of NLS and TL, respectively. These models are used by the optimization unit to compute the objective function and by the logging module of the simulation unit. The logged data are later used to assess the effectiveness of MUTE-DSS, as shown in Fig. \ref{fig_pipeline}.

\subsubsection{SRKW modeling}\label{sec_srkw_modeling}
The SRKW model is broken down into two main tasks: initializing the mammal nodes and assigning trajectories to each. The MUTE-DSS allows the operator to specify the desired number of mammals in the 3D environment, while the modeling unit determines the exact locations and trajectories. As described in Section \ref{sec_database}, the available sighting data and the maximum diving depth data are used to fit a kernel density estimator (KDE), which is then used to initialize the mammal nodes independently. This approach enables data-informed positioning of each mammal.

Let $(\Omega, \mathcal{F}, P)$ be a probability space and let $\mathbf{X}_i : \Omega \rightarrow \mathbb{R}^d$ for $i = 1, 2, \ldots, n$, where $d \geq 1$, be a sequence of independent and identically distributed random variables with an unknown density function $f$. Given mammal sighting data and depth estimates $\{\mathbf{x}_i\}_{i=1}^n$ (realizations of $\{\mathbf{X}_i\}_{i=1}^n$), our aim is to estimate $f$ \cite{silverman2018density}.
Let $\{h(n)\}_{n=1}^\infty \subset (0, +\infty)$ with $h(n) \rightarrow 0$ as $n \rightarrow \infty$. Let $\phi : \mathbb{R}^d \rightarrow \mathbb{R}$ be a measurable, nonnegative kernel function. The KDE is defined by:
\begin{equation}
\mathrm{KDE}_h(\mathbf{x}) = \frac{1}{n h^d(n)} \sum_{i=1}^n \phi\left( \frac{\mathbf{x} - \mathbf{x}_i}{h(n)} \right), \quad \mathbf{x} \in \mathbb{R}^d,
\end{equation}
where $h$ is the bandwidth (also known as the width of the window or the smoothing parameter). The kernel function $\phi$ influences the estimator's overall shape, regularity, and symmetry to approximately zero, while the bandwidth $h$ determines the amount of smoothing. In particular, $\mathrm{KDE}_h$ is a density function provided that $\phi \geq 0$ and $\int_{\mathbb{R}^d} \phi(\mathbf{t}) \, d\mathbf{t} = 1$. Fig.~\ref{fig_SRKW} shows the KDE used to fit the geospatial distribution data and depth distribution data.

The factors driving animal movement and resource selection are complex and diverse. Numerous movement models exist, with a primary differentiation being whether the movement process is framed in continuous or discrete time. The contrasts and similarities between these continuous and discrete-time mechanistic movement models are presented in \cite{mcclintock2014discrete}. This thesis developed a predictive movement model for SRKW aimed at providing timely warnings to reduce collision risks with commercial vessels and limit noise exposure \cite{lin2025forecasting}. Using archival sightings from the Orca Master dataset, which contains reasonably accurate historical movement paths for SRKWs, a velocity-driven Ornstein-Uhlenbeck process was used to predict the mammal's position over time given the initial position. A coastline avoidance algorithm and a direction-blend algorithm were incorporated into the O-U process. Real-time animal movement forecasting has advanced with the integration of state-space models and ensemble-based sequential Monte Carlo methods, enabling probabilistic trajectory predictions. A recent study have demonstrated the effectiveness of this framework in predicting the movements of SRKWs up to 2.5 hours and spatial accuracy of 5 km \cite{randon2022real}. These models, however, is not incorporated into the MUTE-DSS. Instead, the trajectory is assigned as a random directional vector with a reasonable velocity based on typical SRKW movement.

\subsection{Optimization unit}\label{sec_optimization_unit}
By establishing the models of ship noise generation, acoustic propagation, and mammal behavior, the system is ready to perform informed voyage planning. The optimization unit operationalizes this by computing the optimal route and speed profile to minimize URN exposure.
The purpose of the optimization unit is to utilize the aforementioned modeling units and provide the optimal operating conditions of the ship to the ship node. This unit consists of two sub-units: the route planning unit and the speed optimization unit, as shown in Fig. \ref{fig_pipeline}. These two units work hand in hand to guide the ship node.
Consider a case where the ship node is at the departure coordinate and needs to reach a destination coordinate within the ETA (hrs). It sends a message to the optimization unit, requesting the route and the respective vessel speeds along that route. The optimization unit begins with the route planning unit, which is described in the next section.

\subsubsection{Route planning unit}\label{sec_route_planning_unit}

Route planning is implemented using the Open Motion Planning Library (OMPL), an open source library for sampling-based motion planning. It offers various sampling-based algorithms \cite{sucanOMPL}. We define the state space $X$, which encompasses all possible configurations, and we use a planning algorithm to find the optimal route defined by a series of states. In MUTE-DSS, SE2StateSpace is used, consisting of a 2D Euclidean space representing the ship's position and a yaw angle representing its orientation. The geographic coordinates of the ship are transformed into this space based on a reference point (departure coordinates). Due to the presence of islands, which are considered obstacles, certain configurations are invalid and represented by the subset $X_{\text{obs}} \subset X$. The set of valid and collision-free configurations is represented by $X_{\text{free}}$.  A state validity checker is implemented based on islands and shipping lanes in the study area to ensure that ship states are sampled in the appropriate region.

For the ship node, we designate a starting position $x_{\text{start}} \in X_{\text{free}}$ and a target region $X_{\text{goal}} \subset X_{\text{free}}$. A route $\gamma : [0, 1] \to X_{\text{free}}$ is defined as a parametric function where the parameter $\tau \in [0,1]$ represents the progression along the route, such that:
\begin{equation}
\begin{aligned}
    \gamma(0) = x_{\text{start}}, \quad \gamma(1) \in X_{\text{goal}}, \quad \gamma(\tau) &\in X_{\text{free}}, \\
     \forall \tau &\in [0, 1]
\end{aligned}
\end{equation}

The route objective function $\mathcal{J}_r : \Sigma_{X_{\text{free}}} \to \mathbb{R}_{\geq 0}$ assigns a non-negative cost to each route. For computational purposes, we discretize the continuous route into $N_w$ waypoints (ship states) and approximate the total cost as:
\begin{equation}
\mathcal{J}_r(\gamma) \approx \sum_{i=1}^{N_w} \mathcal{J}_{\text{local}}(\gamma(\tau_i)) \Delta \tau_i
\end{equation}
where $\mathcal{J}_{\text{local}}(\gamma(\tau_i))$ is the local cost at waypoint $i$, and $\Delta \tau_i$ is the parameter increment.
The optimization problem is the following.
\begin{equation}
\begin{aligned}
\min_{\gamma} \quad & \mathcal{J}_r(\gamma) \\
\text{Subject to} \quad &
\begin{cases}
\gamma(0) = x_{\text{start}}, \\
\gamma(1) \in X_{\text{goal}}, \\
\gamma(\tau) \in X_{\text{free}}, \quad \forall s \in [0,1].
\end{cases}
\label{eq_route_opt}
\end{aligned}
\end{equation}
In MUTE-DSS the local cost function is defined as the negative mean of transmission losses to the $M$ marine mammal nodes in states ${m}_j \in X_{\text{free}}$:
\begin{equation}
\mathcal{J}_{\text{local}}(\gamma(\tau_i)) = - \frac{1}{M} \sum_{j=1}^{M} \text{TL}(\gamma(\tau_i), {m}_j)
\end{equation}

The optimization problem in Eq. \ref{eq_route_opt} is solved using the Batch Informed Trees (BIT*) algorithm \cite{gammell2015batch}. BIT* performs an ordered search on a continuous planning domain by processing batches of uniformly distributed samples from $X_{\text{free}}$. The algorithm constructs an implicit random geometric graph with samples $X_{\text{samples}} \subset X_{\text{free}}$ and builds an explicit tree outward from $x_{\text{start}}$ towards $X_{\text{goal}}$ using admissible heuristics $\hat{g}(x)$ and $\hat{h}(x)$ to estimate cost-to-come and cost-to-go, respectively. By pruning suboptimal sections of the tree after each batch, BIT* keeps its data structures lean, avoiding unnecessary expansions. Furthermore, drawing samples in batches allows the algorithm to adaptively refine the search region based on the current best solution, reducing the planning space. These features make BIT* converge more quickly to high-quality solutions compared to the previous methods.
The route planning unit provides a series of $N_w$ optimal ship waypoints covering a total distance of $l$ (NM), which must be covered by a speed profile that minimizes the impact of the URN, as explained in the following section.

\subsubsection{Speed optimization unit}\label{sec_speed_optimization_unit}

The optimal ship's route is segmented into $N_t$ sailing legs of fixed duration $\Delta t = \text{ETA} / N_t$ (hrs). Each sailing leg has a specific ship speed $v_i$ (kt), which is the decision variable. Given the speed and duration of each sailing leg, the total distance traveled (TDT) in NM by the ship is defined as:
\begin{equation}
\text{TDT} = \sum_{i=1}^{N_t - 1} v_i \Delta t
\end{equation}

For a single mammal node, the peak received noise level $\text{NL}_i$ is determined by the difference between the noise level at the source $\text{NLS}_i$ and the transmission loss $\text{TL}_i$ of the noise signal as described in Section \ref{sec_acoustic_modeling}. The overall or the broadband received noise level in dB re $1~\mu$Pa for each sailing leg is calculated as:
\begin{equation}
\text{NL}_i(v_i) = \text{NLS}(v_i) - \text{TL}(v_i)
\end{equation}
The root mean square sound pressure $p_i$ for each leg is computed from $\text{NL}_i$ using:
\begin{equation}
p_i = p_{\text{ref}} \times 10^{\frac{\text{NL}_i}{20}}
\end{equation}
where $p_{\text{ref}} = 1~\mu\text{Pa}$ is the reference sound pressure underwater.
The sound exposure $E_i$ for each sailing leg $i$ is calculated by:
\begin{equation}
\begin{aligned}
E_i &= p_i^2 \Delta t = \left( p_{\text{ref}} \times 10^{\frac{\text{NL}_i}{20}} \right)^2 \Delta t = p_{\text{ref}}^2 \times 10^{\frac{\text{NL}_i}{10}} \Delta t
\end{aligned}
\end{equation}
The total sound exposure $E_{\text{total}}$ for the entire voyage is the sum of the exposures from all sailing legs:
\begin{equation}
E_{\text{total}} = \sum_{i=1}^{N_t - 1} E_i = p_{\text{ref}}^2 \Delta t \sum_{i=1}^{N_t - 1} 10^{\frac{\text{NL}_i}{10}}
\end{equation}
The total sound exposure level $\text{SEL}_{\text{total}}$ is defined as:
\begin{equation}
\text{SEL}_{\text{total}} = 10 \log_{10} \left( \frac{E_{\text{total}}}{E_{\text{ref}}} \right)
\end{equation}
where $E_{\text{ref}} = p_{\text{ref}}^2 \times t_{\text{ref}}$ and $t_{\text{ref}} = 1~\text{hr}$ are the reference energy and time duration, respectively.
Substituting $E_{\text{total}}$ and $E_{\text{ref}}$, the total sound exposure level experienced by a single mammal for the entire voyage is given by:
\begin{equation}
\begin{aligned}
\text{SEL}_{\text{total}} &= 10 \log_{10} \left( \frac{p_{\text{ref}}^2 \Delta t \sum_{i=1}^{N_t - 1} 10^{\frac{\text{NL}_i}{10}}}{p_{\text{ref}}^2 t_{\text{ref}}} \right) \\
&= 10 \log_{10} \left( \frac{\Delta t}{t_{\text{ref}}} \sum_{i=1}^{N_t - 1} 10^{\frac{\text{NL}_i}{10}} \right)
\end{aligned}
\end{equation}

The speed optimization problem constituted by the aforementioned SEL averaged over all the mammals and inequality constraints is formulated as:
\begin{equation}
\begin{aligned}
\min_{\mathbf{v} \in \mathbb{R}^{N_t - 1}} \quad & \mathcal{J}_s(\mathbf{v}) = \frac{1}{M} \sum_{j=1}^{M} \text{SEL}_{\text{total},j}(\mathbf{v}) \\
\text{Subject to} \quad &
\begin{cases}
|\text{TDT} - l| \leq \epsilon, \\
v_i-v_{\max } \leq 0, \quad i=1, \cdots, N_t-1 \\
v_{\min }-v_i \leq 0, \quad i=1, \cdots, N_t-1  
\end{cases}
\end{aligned}
\end{equation}
where \( \mathcal{J}_s \) is the speed objective function that represents the mean SEL experienced by \( M \) SRKWs in dB re \( 1 \mu\text{Pa}^2 \text{h} \), TDT is the total distance traveled by ship in NM, \( l \) is the length of the path from the route planning unit in NM, \( v_i \) is the vessel speed at the \( i^{\text{th}} \) sailing leg (kt), and \( v_{\min} \) and \( v_{\max} \) are the minimum and maximum ship speeds (kt), respectively, and finally, $\epsilon = 100$ m is the tolerance for the distance constraint.

The genetic algorithm (GA) is implemented in MUTE-DSS using the Pymoo Python library \cite{blank2020pymoo}. The GA parameters include the number of generations, the population size, and adaptive crossover and mutation coefficients. A population size of 1000 is used, while the crossover and mutation rates dynamically adjust throughout the iterations. This adaptive strategy promotes broad exploration in the early stages and a more focused search toward convergence. The termination criterion is based on either an objective space tolerance or a predefined time window specific to each case study. For the voyage constraints, a niching penalty approach \cite{dobnikar1999niched} is used to enforce the TDT constraint.
The results of the optimization process are used to guide the movement of the ship and simulate interactions with marine mammals. These operations are handled by the simulation and visualization unit, which forms the interactive core of MUTE-DSS.

\subsection{Simulation and visualization unit}\label{sec_simulation_visualization_unit}
The simulation framework is a ROS2-centric system \cite{macenskirobot} where the ship and mammals are represented as individual nodes, as illustrated in Fig. \ref{fig_pipeline}. These nodes are programmed to exhibit specific behaviors within a 3D environment. The ship node receives navigation commands from the optimization unit on the route to follow and the speed to maintain along that route as it travels from its current location to the destination. On the other hand, mammal nodes are launched in the 3D environment and receive behavior directives from the SRKW modeling unit. The states of these nodes are continuously utilized by the optimization unit to dynamically update its solution. This creates a closed-loop system with iterative interactions between nodes and optimization units, synchronized by the time scheduler. The simulation unit also exports the states of the ship and mammal nodes with respect to the simulated voyage time for post-processing.

The visualization unit leverages RViz2 \cite{macenskirobot} for 3D rendering, providing real-time visualization of the noise level contours relative to the ship's location. In addition, time-series data, such as vessel speed and NLS, are plotted in real-time. The visualization unit is interactive, enabling functionalities like zooming into regions of interest, displaying the ship's velocity vector, and overlaying the planned route within RViz2. 
With the core components of the MUTE-DSS framework established, we now assess its performance through a set of realistic case studies based on actual vessel traffic and regional ecological data in the Salish Sea.

\section{Case studies}\label{sec_case_studies}

\begin{table*}[ht]
    \centering
    \resizebox{0.8\textwidth}{!}{%
    \begin{tabular}{l|cccc|cc}
    \hline
    \multicolumn{1}{c|}{Case}  & \multicolumn{4}{c|}{Ship characteristics} &
    \multicolumn{2}{c}{Voyage configurations} \\ 
    \cline{2-7}
    \multicolumn{1}{c|}{study} & Name & Type & \begin{tabular}[c]{@{}c@{}}Length\\(ft)\end{tabular} & 
    \begin{tabular}[c]{@{}c@{}}Speed limits\\(kt)\end{tabular} & 
    \begin{tabular}[c]{@{}c@{}}TDT \\(NM)\end{tabular} & 
    \begin{tabular}[c]{@{}c@{}}ETA\\(h)\end{tabular} \\ \hline

    C1 & \begin{tabular}[c]{@{}c@{}}Star\\Kirkenes\end{tabular} & Other & 684.97 & 8.0 - 16.0 & 156.91 & 12.36 \\
    C2 & Okiana & Other & 697.19 & 8.0 - 14.5 & 124.00 & 9.76 \\ 
    \hline
    \end{tabular}%
    }
    \caption{Summary of ship characteristics and voyage configurations used in the case studies.}
    \label{tab_cases}
\end{table*}

\begin{table*}[ht!]
    \centering
    \resizebox{\textwidth}{!}{%
    \begin{tabular}{cccccccc}
    \hline
    \multicolumn{1}{l}{\multirow{2}{*}{Case study}} 
    & \multicolumn{3}{c}{AIS}                                                                                                        
    & \multicolumn{3}{c}{Optimal}                                                                                                     
    & \multirow{2}{*}{\begin{tabular}[c]{@{}c@{}}$\Delta$ $ \mathcal{J}_s$ \\ (\(\text{dB ref 1} \mu\text{Pa}^2 \text{h} \))\end{tabular}} \\
    \multicolumn{1}{l}{}                            
    & \begin{tabular}[c]{@{}c@{}}ETA \\ (h)\end{tabular} 
    & \begin{tabular}[c]{@{}c@{}}TDT \\ (NM)\end{tabular} 
    & \begin{tabular}[c]{@{}c@{}} $ \mathcal{J}_s$ \\ (\(\text{dB ref 1} \mu\text{Pa}^2 \text{h} \))\end{tabular} 
    & \begin{tabular}[c]{@{}c@{}}ETA \\ (h)\end{tabular} 
    & \begin{tabular}[c]{@{}c@{}}TDT\\  (NM)\end{tabular} 
    & \begin{tabular}[c]{@{}c@{}} $ \mathcal{J}_s$ \\ (\( \text{dB ref 1} \mu\text{Pa}^2 \text{h} \))\end{tabular} 
    &                                                                                    \\ \hline
    C1M1 & 12.36  & 156.91  & 117.81  & 12.36  & 155.76  & 110.67  & -7.14 \\
    C1M2 & 12.36  & 156.91  & 118.55  & 12.37  & 155.14  & 113.65  & -4.90 \\
    C2M1 & 9.76   & 124.00  & 114.85 & 9.75    & 122.16   & 110.34  & -4.51    \\
    C2M2 & 9.76   & 124.00  & 116.40  & 9.75    & 122.20   & 115.48  & -0.92     \\ \hline
    \end{tabular}%
}
    \caption{Comparative analysis of voyage constraints and noise exposure levels for all case studies.}
   \label{tab_results}
\end{table*}

This section presents case studies conducted to demonstrate the applicability of MUTE-DSS and its comparison with real-world ship voyages. A total of two case studies are analyzed, benchmarking their results against the MUTE-DSS optimal voyage. The selected case studies are (i) Star Kirkenes and (ii) Okiana, all departing from the Port of Vancouver, Canada.

\begin{figure}[h!]
\begin{center}
\includegraphics[width=0.45\textwidth]{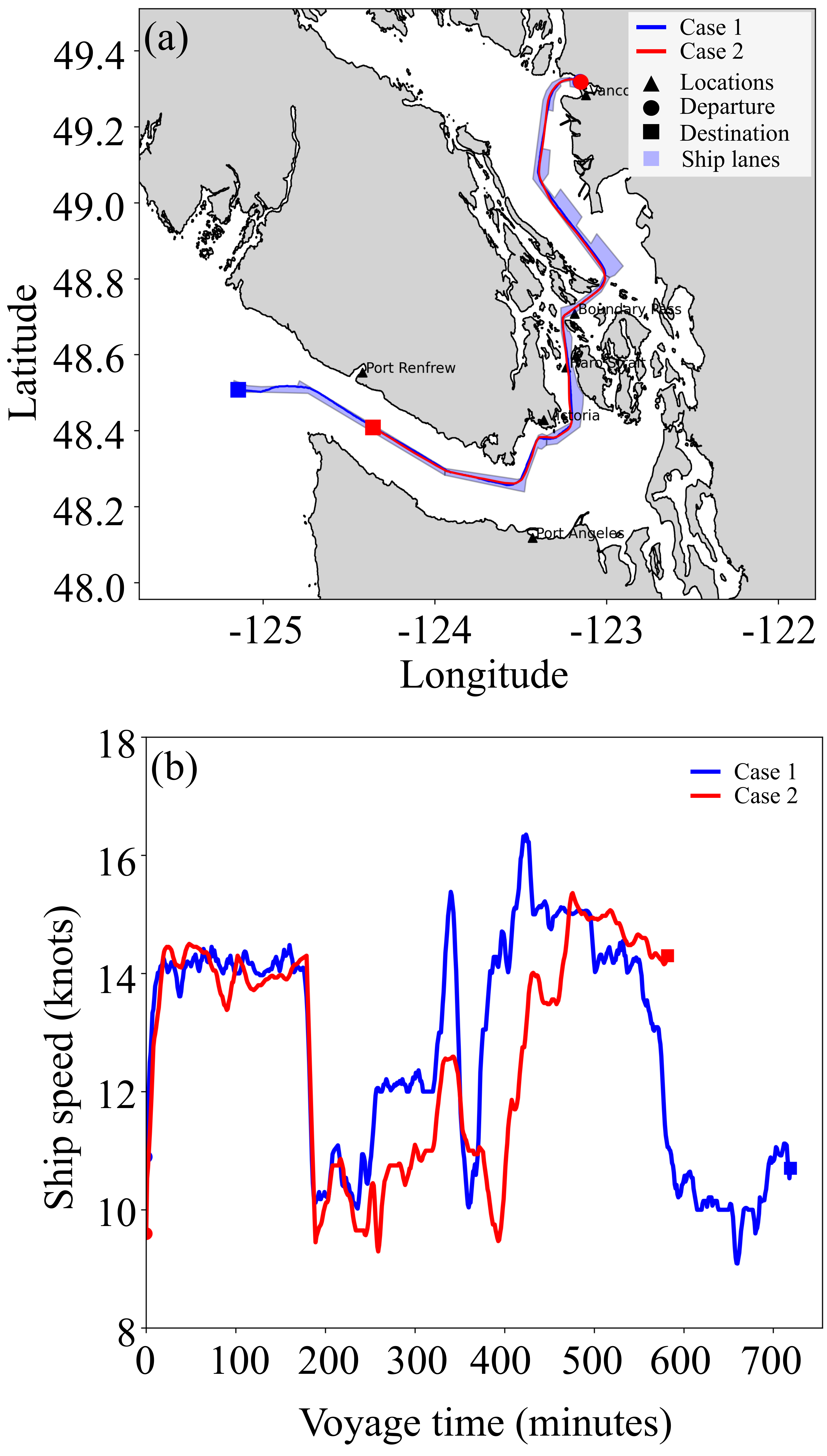}
    \caption{Case study details: (a) ship routes from AIS-GPS; (b) ship speed profile over voyage time from AIS-SOG.}
    \label{fig_casestudy}
\end{center}
\end{figure}

To simulate these voyages within the MUTE-DSS framework, AIS data was processed to extract the route and speed profiles. AIS data includes the following information: ship positions derived from GPS, course, heading, rotation angle, speed, loading status, ship type, navigation status, destination, and timestamp. Missing entries in the raw data were interpolated, and the data were reformatted to ensure consistency. Geographic coordinates and vessel speed were tracked relative to voyage time. The extracted routes and speed profiles for the two ships are depicted in Fig. \ref{fig_casestudy}, while each ship's characteristics are summarized in Table \ref{tab_cases}. The voyage constraints of the route and speed optimization problems are also derived from these data for each case study.
The case studies described above serve as the basis for quantitative analysis. In the following section, we evaluate the acoustic impact of AIS-based versus optimized voyages and interpret the system’s ability to adapt to diverse ecological configurations.

The case studies were simulated under two distinct configurations of SRKW presence within a 3D environment:
\begin{enumerate}
    \item M1: A single mammal node was positioned at latitude 48.646343\degree N, longitude 123.313054\degree W, and a depth of 1 meter. By constraining the mammal’s position, we isolate and assess the MUTE-DSS’s capacity to solve the voyage optimization problem for a single, stationary receptor node.
    \item M2: Five mammal nodes were modeled as described in Section \ref{sec_srkw_modeling}. These simulated nodes move dynamically with respect to the simulated voyage time. This configuration represents a more complex and realistic real-world scenario.
\end{enumerate}

\section{Results and discussion}
\label{sec_results}


This section presents the results of case studies conducted under two SRKW configurations (M1 and M2). For each configuration, simulations were run for both AIS‑based voyages and voyages with optimized operating conditions. The results include a comparative assessment of the route trajectories and speed profiles of the ships. In addition, an acoustic footprint analysis quantifies the levels of noise exposure at the SRKW nodes and compares the exposure between the two types of voyage. Table \ref{tab_results} summarizes these findings.

\begin{figure*}[ht]
\centering
\includegraphics[width=\textwidth]{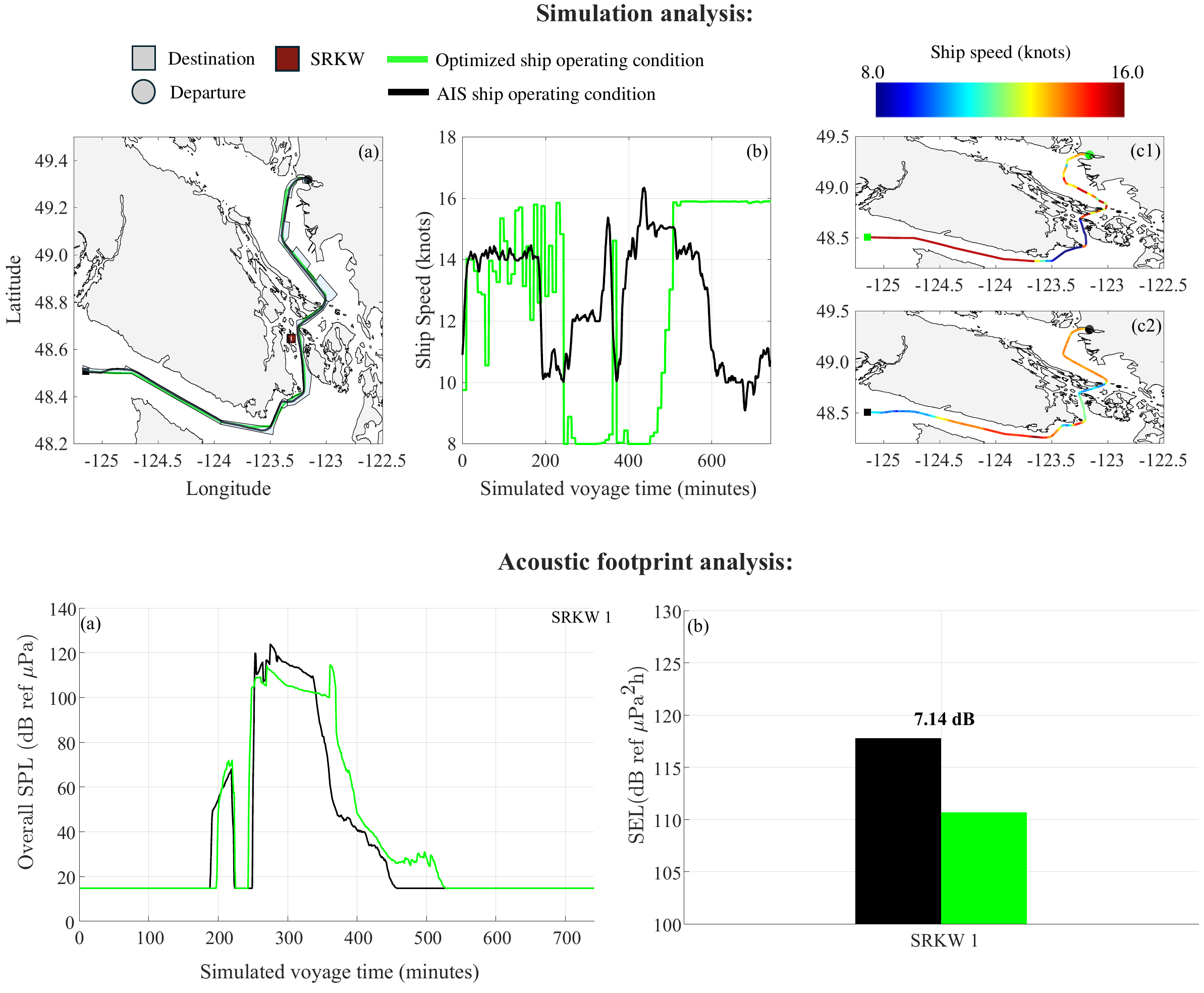}
\caption{Results for Case 1 in the M1 configuration. Voyage comparison between AIS-based and optimized paths: (a) route, (b) speed profiles, and (c) speed-color-coded routes. Acoustic footprint analysis: time series of (a) overall SPL and (b) noise exposure levels.}
\label{fig_ch6_C1M1}
\end{figure*}

Figure \ref{fig_ch6_C1M1} shows the results for Case C1 in the SRKW configuration M1. In subplot (a), the route derived from AIS data is plotted along the route optimized by MUTE-DSS; both routes remain within designated shipping lanes and cover distances of 156.91 NM and 155.76 NM, respectively. Notably, the optimized route departs from the AIS track at several segments, which reflects the route planning unit’s goal of increasing transmission loss to the SRKW nodes while respecting voyage constraints. Since TL increases with distance when environmental factors are fixed, the planner deliberately maximizes the separation between vessel and mammal to improve noise attenuation.

Subplot (b) of Fig. \ref{fig_ch6_C1M1} compares ship speed throughout the simulated voyage. In addition, subplots (c1) and (c2) map the optimized route with speed represented by a color gradient. Initially, up to roughly 250 minutes, the ship speed fluctuates between 11 kt and 16 kt; this variability arises from the interaction of URN and bathymetric features. In this segment, the presence of an archipelago (the Galiano islands) between the ship and mammal nodes creates alternating high‑intensity and shadow zones. Although the Gaussian RBF interpolant used in the optimization unit captures TL at an acceptable resolution, a higher‑fidelity model might yield a smoother speed profile. At approximately 250 minutes of simulated voyage time, speed abruptly decreases to 8 kt, reflecting the absence of intervening landmasses to attenuate URN. Thereafter, a pronounced spike followed by a sharp deceleration again corresponds directly to the bathymetric features. Finally, beyond 500 minutes, the ship accelerates to 16 kt upon entering the Juan de Fuca Strait, where Vancouver Island fully obstructs noise propagation from the SRKW nodes. These characteristics in the speed profile indicate that MUTE-DSS leverages the shadow zone to make informed decisions about adjusting the ship’s speed.

The acoustic footprint analysis is also presented in Fig. \ref{fig_ch6_C1M1}. During each simulation, received noise levels at the mammal nodes are continuously recorded as SPL in dB re 1 $\mu$Pa. The left‑hand subplot (a) shows the SPL time series for both the AIS‑based and optimized voyages. The maximum SPL during the AIS‑based voyage is considerably higher than during the optimized voyage. However, the optimized voyage sustains lower SPL levels for a longer duration, producing a temporally stretched time series compared to the other. The right‑hand subplot (b) depicts the corresponding sound exposure levels (SEL in dB re 1 $ \mu\text{Pa}^2$ h). Notably, the optimized route reduces overall exposure by 7.14 dB, equivalent to an approximate five‑fold decrease, resulting in only 19.32\% of the exposure observed in the AIS‑based smulation. These results demonstrate the effectiveness of MUTE‑DSS in mitigating marine mammals’ exposure to URN.

\begin{figure*}[ht]
\centering
\includegraphics[width=\textwidth]{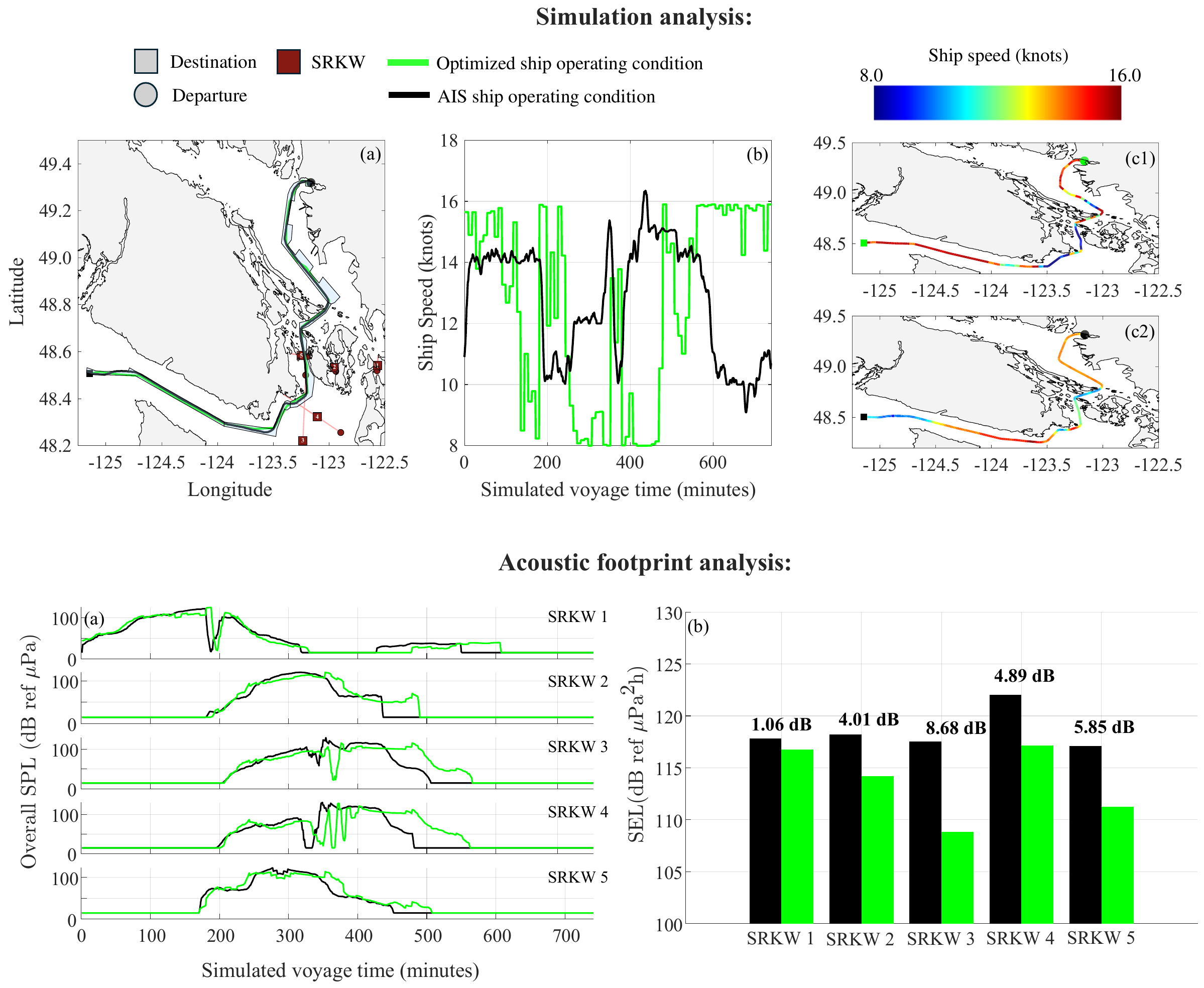}
\caption{Results for Case 1 in the M2 configuration. Voyage comparison between AIS-based and optimized paths: (a) route, (b) speed profiles, and (c) speed-color-coded routes. Acoustic footprint analysis: time series of (a) overall SPL and (b) noise exposure levels.}
\label{fig_ch6_C1M2}
\end{figure*}

Furthermore, a supplementary case study examined the separate contributions of ship route and vessel speed to SEL. In this case, the ship adhered exactly to the AIS route, enabling the optimization unit to adjust only the ship speed. Under these constraints, the MUTE DSS achieved a reduction of 6.47 dB in SEL, compared with the previous reduction of 7.14 dB. Thus, it can be concluded that the route modifications accounted for only 0.67 dB of the overall decrease.

In configuration M2, five mammal nodes are distributed according to the method outlined in Section \ref{sec_srkw_modeling}. Figure \ref{fig_ch6_C1M2} presents both vessel trajectories and noise exposure analysis for this case. In the first geographic subplot (a), the spatial distribution and trajectories of the five mammal nodes are shown. As in the previous case, the optimized route departs only slightly from the AIS route. The optimal speed profile, shown in subplot (b), follows a generally similar trend, with three major slowdowns: the first, occurring before entering Boundary Pass, spans 100 to 200 minutes of simulation time and primarily targets SRKW1. Moreover, similar to the behavior observed in the previous case, the ship speed oscillates while transiting the archipelago in the Strait of Georgia. Subsequently, a second significant slowdown occurs as the ship enters Haro Strait, covering the voyage time from 250 to 500 minutes, with a peak at around 375 minutes. These speed reductions are identified by the algorithm to mitigate exposure, as most SRKWs are in close proximity to the ship in this region. Beyond the 500‑minute mark, there is the third slowdown to 12 kt, specifically targeting SRKW3, as illustrated in subplot (a), before cruising again at a peak speed of 16 kt.

\begin{figure*}[ht]
\centering
\includegraphics[width=\textwidth]{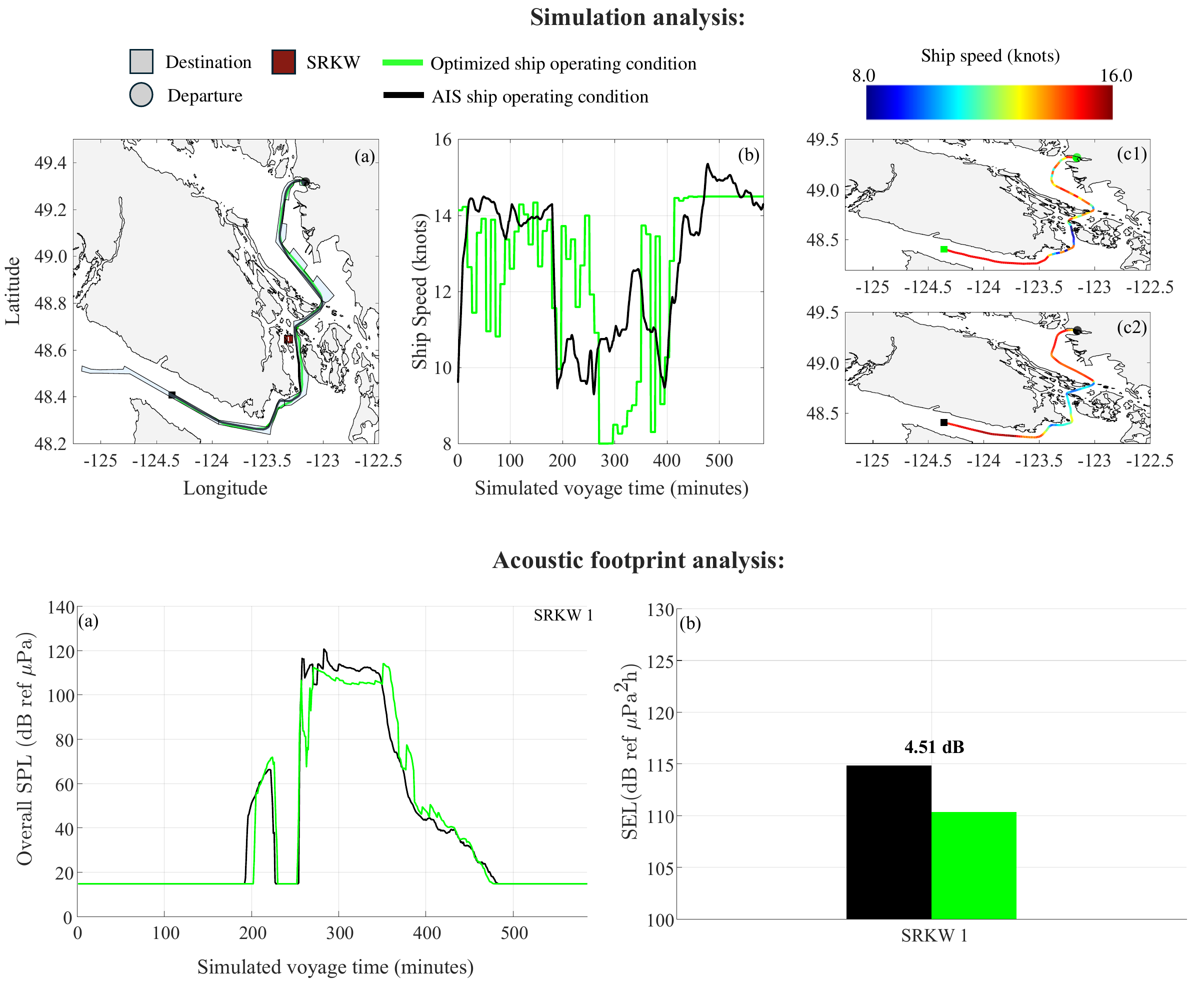}
\caption{Results for Case 2 in the M1 configuration. Voyage comparison between AIS-based and optimized paths: (a) route, (b) speed profiles, and (c) speed-color-coded routes. Acoustic footprint analysis: time series of (a) overall SPL and (b) noise exposure levels.}
\label{fig_ch6_C2M1}
\end{figure*}

The acoustic footprint analysis in Fig. \ref{fig_ch6_C1M2} reveals a consistent reduction in SEL received by all SRKWs throughout the optimized voyage, indicated by the green bar. On average, the AIS-based voyage exhibits an SEL of 118.55 dB, whereas the optimized voyage achieves an SEL of 113.65 dB, resulting in a reduction of 4.90 dB, or roughly only 32.4\% of the AIS-based SEL. Subplot (b) of the acoustic analysis highlights that some SRKWs experience greater reductions than others. Specifically, the third SRKW shows the most reduction, with SEL decreasing by 8.68 dB. This reduction is attributed to its trajectory within the environment, as depicted in the geo subplot (a), where this SRKW is initialized near the ship lane but moves away. In contrast, the first SRKW exhibits the least reduction, with a decrease of only 1.06 dB compared to the AIS-based voyage. This can be explained by its distant initial position relative to the voyage, as shown in the geo subplot (a). Its prolonged exposure to open waters during the entire voyage results in nearly equivalent SEL contributions between the optimized and AIS-based scenarios. Therefore, MUTE-DSS effectively reduces the acoustic footprint by leveraging the TL model, which identifies shadow zones. However, if a marine mammal is exposed in open water, it is likely to experience levels of exposure comparable to or even higher than in the AIS-based voyage.

Figure \ref{fig_ch6_C2M1} presents the results of Case C2 for the first SRKW configuration, M1. As shown in Table \ref{tab_results}, this case has shorter ETA and TDT compared to the previous case. Consistent with previous studies, the optimized route and speed profiles closely resemble those observed earlier. However, a key difference is that the slowdown interval is reduced to approximately 100 minutes of simulation time, rather than the 200-minute interval noted in prior cases. This reduction stems primarily from the vessel’s speed constraint of under 14.5 kt, which precluded a prolonged slowdown in order to satisfy the ETA requirement. The subplot (c1) shows that this deceleration occurs in close proximity to the SRKW location.

\begin{figure*}[ht]
\centering
\includegraphics[width=\textwidth]{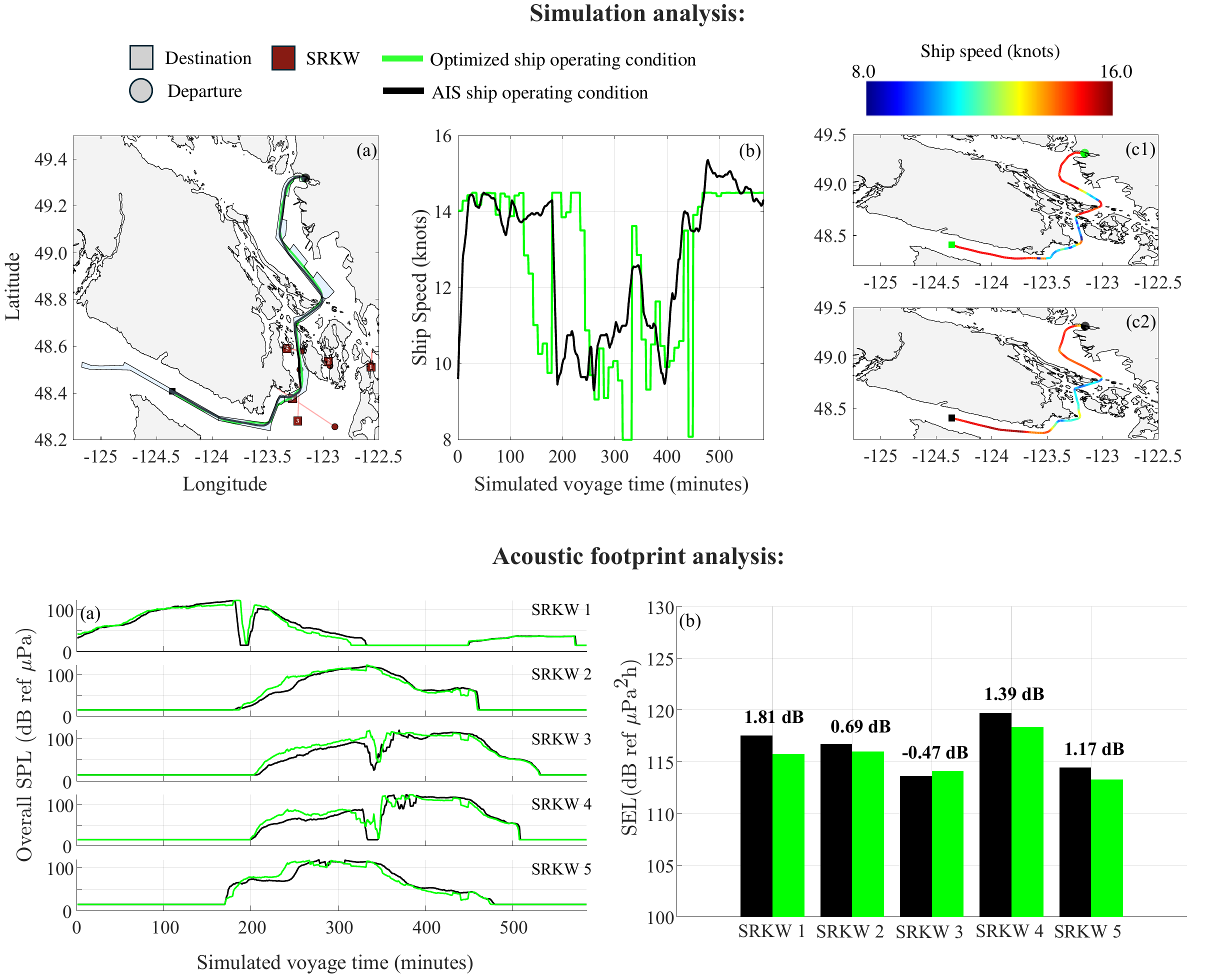}
\caption{Results for Case 2 in the M2 configuration. Voyage comparison between AIS-based and optimized paths: (a) route, (b) speed profiles, and (c) speed-color-coded routes. Acoustic footprint analysis: time series of (a) overall SPL and (b) noise exposure levels.}
\label{fig_ch6_C2M2}
\end{figure*}

The acoustic footprint analysis in Fig. \ref{fig_ch6_C2M1}  further complements this characteristic. Unlike the gradually decreasing SPL observed in the previous case (Fig. \ref{fig_ch6_C2M1}), the overall SPL time series for the optimized voyage maintains relatively consistent levels throughout the exposure duration. This is attributed to the absence of a prolonged slowdown. Nevertheless, the peak SPL remains lower than that of the AIS-based voyage, demonstrating the effectiveness of the algorithm. As a result, the SEL is reduced by 4.51 dB, achieving a 65.63\% reduction. However, when compared to the optimal solutions from the other two cases within the same SRKW configuration, this case exhibits the highest SEL. This increase is primarily due to tighter time constraints, which limit the extent of speed adjustments available for noise mitigation.

Figure \ref{fig_ch6_C2M2} presents the results of Case C2 for the SRKW configuration, M2. The speed profile again reveals three distinct slowdown regions: the first occurring between 150 and 200 minutes of simulation time before the ship enters Boundary Pass, the second spanning approximately 250 to 325 minutes along Haro Strait, and the final slowdown occurring between 350 and 425 minutes. These features are consistent with the previous case. However, a key distinction in this case is that the speed does not decrease to 8 kt, and the duration of the slowdown is reduced due to a shorter time constraint, as previously discussed.

The acoustic footprint analysis of C2 under M2 configuration indicates that the average SEL is reduced by only 0.92 dB as presented in Table \ref{tab_results}. Fig. \ref{fig_ch6_C2M2} illustrates the individual exposure levels of the SRKWs, showing that while the optimized voyage achieves lower SEL for some SRKWs, others experience slightly higher exposure compared to the AIS-based voyage. The first SRKW undergoes the most significant reduction in SEL, with a decrease of 1.81 dB, attributed to the initial slowdown, which is not implemented in the AIS-based voyage. However, as a consequence, the optimization algorithm does not allocate sufficient time to reduce the ship’s speed to 8 kt in the Haro Strait region, unlike in the other case. Consequently, the SEL values for the remaining mammal nodes are elevated compared with the previous case; in particular, SRKW3 experiences the greatest exposure, with an SEL 0.47 dB higher than that recorded during the AIS‑based voyage. Nevertheless, on average, the optimized voyage generated by MUTE-DSS achieved a modest reduction of 19.11\%.
The observed reductions in sound exposure levels demonstrate the system’s capacity to adaptively plan noise-aware voyages. We conclude by summarizing our key findings and outlining future extensions to broaden the capabilities and impact of MUTE-DSS.

\section{Conclusion}\label{sec_conclusion}
We have presented the Mitigating Underwater Noise Transmission and Effects Decision Support System (MUTE-DSS), a ship voyage planning DSS that minimizes cumulative underwater radiated noise exposure to marine mammals. MUTE-DSS is a ROS2-centric digital twin framework that unifies state-of-the-art acoustic ship noise signature modeling, data-informed SRKW behavior, and two-stage voyage optimization. The acoustic signature model integrated JOMOPAN-ECHO for the ship noise source spectrum with propagation losses that were derived from a 3D Bellhop and modeled using Gaussian radial basis functions. Marine mammals were initialized using kernel density estimates of Southern resident killer whale sightings and dive depths, then advanced with randomized velocity vectors. These were used by a two-stage optimization unit: Batch-Informed Trees generated an optimal collision-free route, and a genetic algorithm provided the optimal ship speed under voyage constraints, minimizing the impact of URN on marine mammals.

The study area covered the region from the Strait of Georgia to the Strait of Juan de Fuca. In this study, two ship voyages (C1 and C2) were simulated with two different marine mammal configurations (M1 and M2) to demonstrate the efficacy of the proposed DSS. The M1 configuration involved a single SRKW at a fixed location at Boundary Pass, while the M2 configuration considered five SRKWs dynamically moving within the environment. The AIS data from a ship operating in the study area were used to construct the baseline AIS-based voyage, which was compared against the optimal voyage generated by MUTE-DSS. The ship routes, speed profiles, and resulting acoustic footprints of the respective voyages were analyzed and compared. 

The results of the simplified configuration, M1, in both cases (C1 and C2) highlight the informed decision-making capabilities of MUTE-DSS by effectively identifying URN propagation patterns within the study area. The study area is characterized by its shallow water environment, where the acoustic propagation is heavily influenced by the structure of the sea floor forming shadow zones. This case illustrates how MUTE-DSS optimization algorithms leverage the shadow zones created by the archipelago to optimize the speed profile. The optimized speed profile involves higher speeds in the shadow zones and reduced speeds when the SRKW is exposed to URN. This significantly reduces the peak SPL received by the SRKW, as observed in the time-series analysis, compared to the AIS-based voyage. Even though the SRKW experienced longer exposure to SPL in the optimal voyage, the substantial reduction in SPL levels resulted in a significant decrease in SEL. Specifically, the cases C1M1 and C2M1 showed SEL reductions of 7.14 dB and 4.51 dB, respectively, corresponding to reductions of 80.68\% and 65.63\% compared to the AIS-based voyages. Furthermore, a special case study, which constrained the ship to the AIS-based route while optimizing only the speed profile, demonstrated that the majority of the SEL reduction can be attributed to speed manipulation.

The results of the realistic configuration, M2, across all cases demonstrate the capability of the MUTE-DSS to account for multiple moving SRKWs. The optimized speed profile has variations throughout the voyage, reflecting the dynamic nature of the SRKWs' movements. As some SRKWs initially in the shadow zones may later become exposed to URN, necessitating adjustments in the ship's speed profile. In the case of C1M2, the SRKW1 experienced an SEL similar to the AIS-based voyage, with a difference of only 1.06 dB. However, the SRKW3 benefited from a substantial reduction of 8.68 dB in SEL. These results highlighted that the efficacy of MUTE-DSS varied across different individuals, depending on their locations and movement patterns. Overall, MUTE-DSS reduced the average SEL by 4.90 dB, corresponding to a 67.6\% reduction. Case C2M2 had tighter constraints than C1M2, including more restrictive speed limits. As a result, the optimal solution produced an acoustic footprint similar to the AIS-based voyage, with only a minor reduction of 0.92 dB on the overall SEL. While MUTE-DSS successfully reduced SEL for some SRKWs, this led to a slight increase in SEL of 0.47 dB for SRKW3, which was exposed toward the end of the voyage. However, these case studies demonstrated the adaptive operational strategy of MUTE-DSS to mitigate the impact of the URN, demonstrating its ability to dynamically adjust ship speed and route while considering environmental and operational factors.

The Integrated Bridge System (IBS) serves as the centralized platform for ship navigation, control, decision support, status monitoring, display, and communication by interconnecting onboard digital devices and acquiring data that describes the vessel’s state. Similarly, to represent the state of marine mammals in the operational area, an external passive acoustic monitoring system can be deployed to detect and localize animals and stream their real-time positions into the IBS. In future work, MUTE-DSS can be incorporated as a module within the IBS, leveraging both ship and marine mammal state inputs to perform voyage planning and formulate navigation strategies that mitigate underwater noise exposure. Using ROS2 as the middleware enables seamless fusion of multimodal sensor data and integration with heterogeneous hardware, ensuring low-latency, extensible coordination between sensing and decision-making.

\section*{Acknowledgments}
The present study is supported by Mitacs, Transport Canada, and Clear Seas through the Quiet-Vessel Initiative (QVI) program. The authors would like to express their gratitude to Dr. Paul Blomerous and Ms. Tessa Coulthard for their valuable feedback and suggestions. In addition, we gratefully acknowledge the contributions of Dr. David Rosen and Dr. Andrew W Trites, who provided insights on the marine biology aspect of the study. This research was also supported in part through computational resources and services provided by Advanced Research Computing (ARC) at the University of British Columbia and Compute Canada.


\appendix

\section{Acoustic model parameters and configurations} \label{sec_app_acoustic}
\begin{figure}[h!]
\begin{center}
\includegraphics[width=0.45\textwidth]{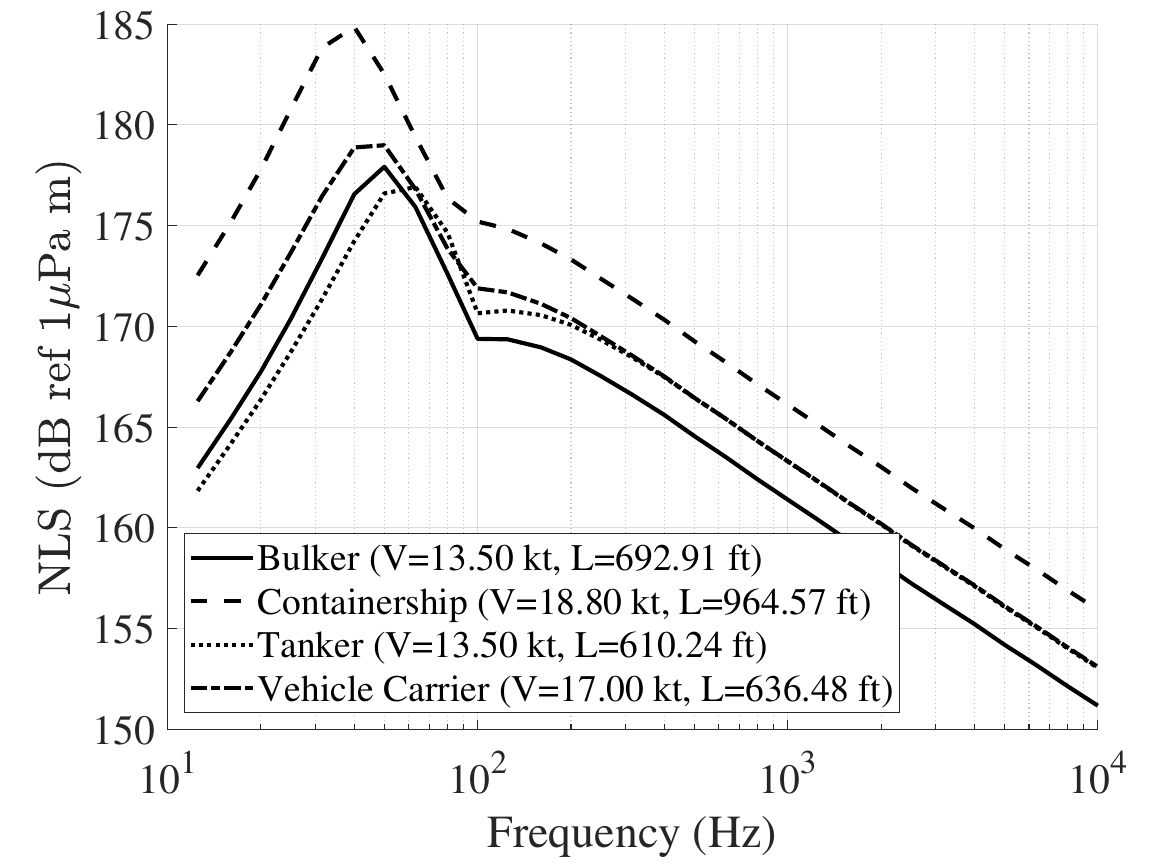}
    \caption{JOMOPANS-ECHO model predictions of decidecade band noise source levels for different ships and operating conditions. }
     \label{fig_jomopan}

\end{center}
\end{figure}

This section details the parametric configuration of the JE noise source and the Bellhop acoustic propagation model implemented in MUTE-DSS. The reference speed ($v_T$) in Eq.~\ref{eq_jomopan} is assigned according to ship type; the values used are summarized in Table~\ref{tab_jomopan}. Figure~\ref{fig_jomopan} presents the resulting NLS spectrum for different combinations of ship type, ship speed, and ship length.

\begin{figure}[h!]
\begin{center}
\includegraphics[width=0.45\textwidth]{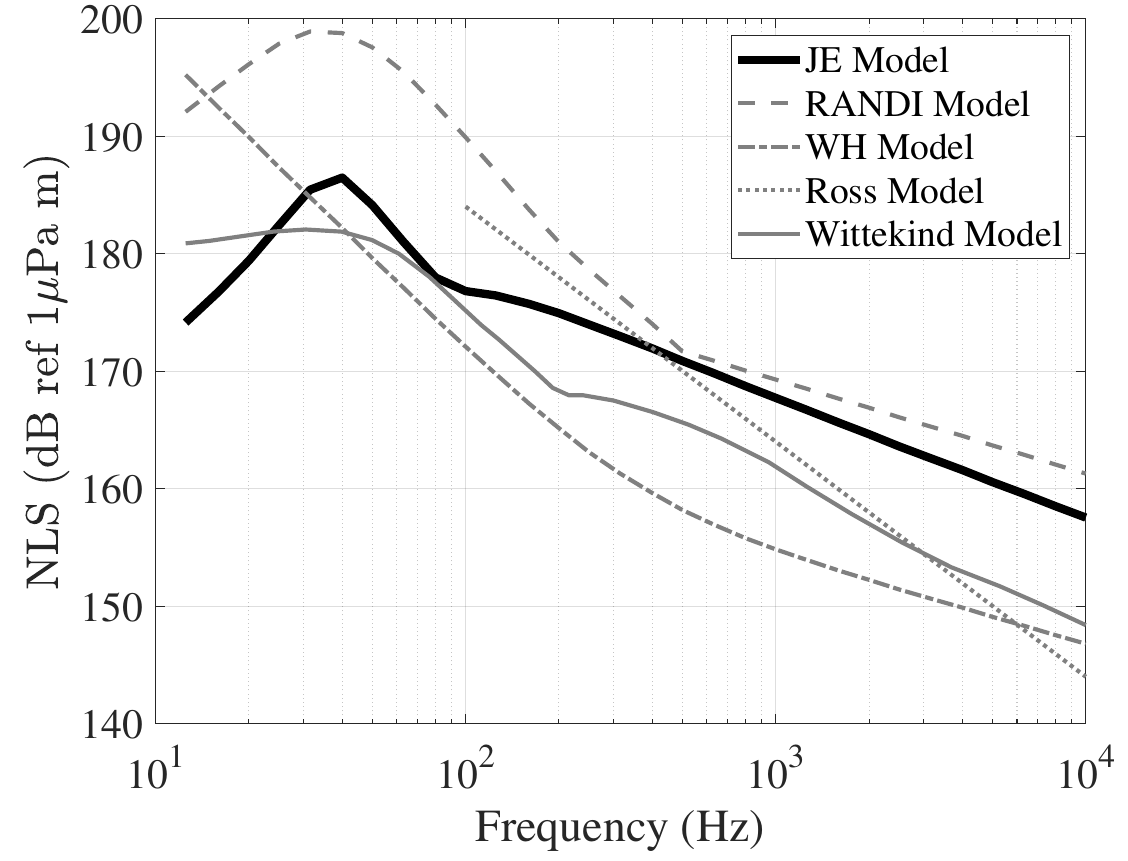}
    \caption{Comparison of ship noise source models: The black solid line represents the JOMOPANS-ECHO model (Eq. \ref{eq_jomopan}) used in the framework. Modeled NLS (dB ref $1\mu$Pa m) are shown for a containership with a displacement of 135,000 tons, a block coefficient of 0.65, a length of 964 ft, and CIS of 14 kt, traveling at 20 kt equipped with two diesel engines, each weighing 50 tons.  }
     \label{fig_NLS_models}
\end{center}
\end{figure}

 Figure \ref{fig_NLS_models} illustrates a comparison of most of the models for a particular ship under specific operating conditions; however, the Urick \cite{urick1983principles} and PIANO \cite{lloyd2024piano} models are not included. Urick's model, being a function of propeller tip speed, could not be modeled in this context. Although the mathematical framework of the PIANO model is provided, the parameters are not shared in the paper. It is evident that Ross's \cite{Ross1976} and Wittekind's models \cite{wittekind2014simple} are the simplest, as they do not incorporate the spectral hump observed in more complex models and have a fixed slope. The RANDI 3.1 model \cite{BreedingJ96}, on the other hand, tends to overestimate noise levels; however, it has a similar high frequency slope to the JE model. Notably, the JE model closely aligns with the Wittekind model at low frequency.

Bellhop is an efficient acoustic ray tracing model implemented in Fortran by Michael Porter \cite{acoustics_toolbox}. It supports both two-dimensional and three-dimensional propagation environments and computes acoustic fields for prescribed sound speed profiles ($c(z)$) or full spatially varying sound speed fields ($c(r,z)$) within waveguides that may have flat or variable boundaries. The model can output a variety of quantities, including ray trajectories, travel times, amplitudes, eigenrays, acoustic pressure, and transmission loss, with coherent, incoherent, or semi-coherent formulations. In the Bellhop environmental configuration used here, the upper ocean surface is assumed to be flat, while the bottom boundary is treated as an acousto-elastic half-space. Geometric Gaussian beams are used \cite{porter1987gaussian}, and coherent transmission loss is computed accordingly. The receiver locations are sampled at up to 2000 positions in both range and bearing, corresponding to a spatial resolution of 50 meters. The beam-launching geometry is discretized with 360 azimuthal ($\beta$) and 360 elevation ($\alpha$) angles, where $\alpha$ spans from $-90^\circ$ to $+90^\circ$ and $\beta$ spans $0^\circ$ to $360^\circ$.

\begin{table}[h!]
  \centering
  \resizebox{\columnwidth}{!}{%
    \begin{tabular}{l l c}
      \hline
      \textbf{Ship Type (T)} & \textbf{AIS Ship Type ID} & \textbf{$v_T$ (kt)} \\ \hline
      Fishing               & 30                       & 6.4                 \\
      Tug                   & 31, 32, 52               & 3.7                 \\
      Naval                 & 35                       & 11.1                \\
      Recreational          & 36, 37                   & 10.6                \\
      Government/Research   & 51, 53, 55               & 8.0                 \\
      Cruise                & 60–69 (length $l > 100$ m)  & 17.1               \\
      Passenger             & 60–69 (length $l \leq 100$ m) & 9.7              \\
      Bulker                & 70, 75–79 ($V \leq 16$ kt)    & 13.9             \\
      Containership         & 71–74, 70, 75–79         & 18.0                \\
      Vehicle Carrier       & n/a                      & 15.8                \\
      Tanker                & 80–89                    & 12.4                \\
      Other                 & All other type IDs       & 7.4                 \\
      Dredger               & 33                       & 9.5                 \\ \hline
    \end{tabular}%
  }
  \caption{Ship types, corresponding AIS IDs, and reference speeds ($v_T$) used in Eq.~\ref{eq_jomopan}.}
  \label{tab_jomopan}
\end{table}

 \bibliographystyle{elsarticle-num-names} 
 \bibliography{biblograph}





\end{document}